\documentclass[12pt]{article}

\usepackage{amsmath}
\usepackage{amsfonts}
\usepackage{amssymb}
\usepackage{graphicx, color}
\usepackage{hyperref}
\usepackage{epsfig}
\usepackage{epstopdf} 

\definecolor{darkgreen}{rgb}{0,0.55,0}

\newtheorem{proposition}{Proposition}[section]
\newtheorem{theorem}{Theorem}[section]
\newtheorem{lemma}[theorem]{Lemma}

\newtheorem{remark}[theorem]{Remark}

\newtheorem{definition}{Definition}

\def\phi{{\varphi}}

\DeclareSymbolFont{AMSb}{U}{msb}{m}{n}
\DeclareMathSymbol{\N}{\mathbin}{AMSb}{"4E}
\DeclareMathSymbol{\Z}{\mathbin}{AMSb}{"5A}
\DeclareMathSymbol{\R}{\mathbin}{AMSb}{"52}
\DeclareMathSymbol{\Q}{\mathbin}{AMSb}{"51}
\DeclareMathSymbol{\I}{\mathbin}{AMSb}{"49}
\DeclareMathSymbol{\C}{\mathbin}{AMSb}{"43}

\DeclareMathOperator*{\esssup}{ess\,sup}
\DeclareMathOperator*{\essinf}{ess\,inf}

\newcommand{\wt}{\widetilde}
\newcommand{\bbl}{[\hspace{-0.07em}[}
\newcommand{\bbr}{]\hspace{-0.07em}]}
\newcommand{\e}{\varepsilon}

\newcommand{\calL}{{\mathcal L}}

\newcommand{\calD}{{\mathcal D}}

\newcommand{\calA}{{\mathcal A}}
\newcommand{\calH}{{\mathcal H}}

\begin{document}

\title{Existence and uniqueness of minimizers of general least gradient problems}

\author{  
{Robert L. Jerrard\footnote{Department of Mathematics, University of Toronto, Toronto, Ontario, Canada M5S
2E4. E-mail: rjerrard@math.toronto.edu.} }\qquad
{Amir Moradifam\footnote{Department of Applied Physics and Applied Mathematics, Columbia University, New York, NY, USA. E-mail: am3937@columbia.edu.}
\qquad Adrian I. Nachman\footnote{Department of Mathematics and the
Edward S. Rogers Sr. Department of Electrical and Computer
Engineering, University of Toronto, Toronto, Ontario, Canada. E-mail: nachman@math.toronto.edu. }\qquad }}

%\thanks{The authors are grateful to Benjamin Stephens for numerous helpful discussions.}

\smallbreak \maketitle

\begin{abstract}
Motivated by problems arising in conductivity imaging, we prove existence, uniqueness, and comparison theorems - under certain sharp conditions -
 for minimizers of the general least gradient problem
\[\inf_{u\in BV_f(\Omega)} \int_{\Omega}\varphi(x,Du),\]
where
$f:\partial \Omega\to \R$ is continuous,
\[
BV_f(\Omega):=\{v\in BV(\Omega): 
\ \ \forall x\in \partial \Omega, \ \ 
\lim_{r\to 0} \ \esssup_{y\in \Omega, |x-y|<r} |f(x) - v(y)| = 0 \ \}
%BV_f(\Omega)=\{u\in BV(\Omega): \hspace{0.1cm} u|_{\partial \Omega}=f \hspace{0.1cm} \hbox{and} \hspace{0.1cm} \hspace{0.1cm} u \hspace{0.1cm} \hbox{is continuous at } \hspace{0.1cm} \partial \Omega \}.
\] 
and $\varphi(x,\xi)$ is a function that, among other properties, is convex and homogeneous of degree $1$ with
respect to the $\xi$ variable.
In particular we prove that if  $a\in C^{1,1}(\Omega)$ is bounded away from zero, then minimizers of the weighted least gradient problem $\inf_{u \in BV_f}\int_{\Omega} a|Du|$ are unique in $BV_f(\Omega)$. We construct counterexamples to show that the regularity assumption $a\in C^{1,1}$ is sharp, in the sense that it can not be replaced by $a\in C^{1,\alpha}(\Omega)$ with any $\alpha<1$.
\end{abstract}

\maketitle

\section{Introduction}
Let $\Omega$ be a bounded open set in $\R^n$ with Lipschitz boundary and $\varphi: \Omega\times \R^n \rightarrow \R$ be a continuous function satisfying the following conditions.\\

($C_1$) There exists $\alpha>0$ such that $\alpha |\xi| \leq \varphi(x,\xi) \leq \alpha^{-1}|\xi|$ for all $x\in \Omega$ and $\xi \in \R^n$. \\

($C_2$) $\xi \mapsto \varphi(x,\xi)$ is a norm for every $x$.\\

For any $u\in BV_{loc}(\R^n)$ let $\varphi(x,Du)$ denote the measure
defined by
\begin{equation}\label{varphi.def}
\int_A \varphi(x,Du) \ = \ \int_A \varphi(x, v^u(x)) |Du|
\qquad\mbox{ for $A$ a bounded Borel set},
\end{equation}
where $|Du|$ is the total variation measure associated to the vector-valued
measure $Du$, and $v^u$ denotes the Radon-Nikodym derivative $v^{u}(x)=\frac{d\, Du}{d\, |Du|}$.
(The right-hand side of \eqref{varphi.def} makes sense, since $v^u$ is $|Du|$-measurable,
and hence $\varphi(x, v^u(x))$ is as well.)
Standard measure-theory considerations and facts about $BV$ functions imply that (see \cite{AB}) if $U$ is an open set,
then
\begin{equation}\label{varphi.def1}
\int_{U} \varphi(x,Du) = \sup \{ \int_{U} u \nabla \cdot Y dx\ \ : \ \ Y \in C^{\infty}_{c} (U; \R^n), \ \ \sup \varphi^0(x, Y(x)) \leq 1 \}, 
\end{equation}
where $\varphi^0(x, \cdot)$ denotes the norm on $\R^n$ dual to $\varphi(x, \cdot)$, defined by
\[
\varphi^0(x,\xi) := \sup \{ \xi \cdot p : \varphi(x,p) \le 1\}.
\]

For $u\in BV(\Omega)$, $\int_{\Omega}\varphi(x, Du)$ is called the $\varphi$-total variation of $u$ in $\Omega$. Also if  $A, E$ are subsets of $\R^n$, with $A$ Borel and $E$ having finite perimeter,
then we shall write $P_{\varphi}(E; A)$ to denote the $\varphi$-perimeter of $E$ in $A$, defined by
\begin{equation}\label{perimeter}
P_{\varphi}(E;A):= \int_A \varphi(x, D \chi_E) ,
%= \sup \{\int_{A} \nabla \cdot Y dx: \ \ Y  \in C^{\infty}_{c} (U; \R^n), \ \ \sup \varphi^0(x, Y(x)) \leq 1\}.
\end{equation}
where $\chi_E$ is the characteristic function of $E$. We will also write 
$P_\varphi(E)$ to mean $P_\varphi(E;\R^n)$.
We remark that if $\partial E$ is smooth enough, then
\[
P_{\varphi}(E;A):= \int_{\partial E \cap A}  \varphi(x, \nu_E(x)) \,d\calH^{n-1}\qquad
\mbox{$\nu_E :=$  outer unit normal,}
\]
which is
a generalized inhomogeneous, anisotropic area of $\partial E$ in $A$.
%A set $E\subset\R^n$ is said to be of locally finite $\varphi$-perimeter  if $P_{\varphi}(E; U)<\infty$ for every bounded open set $U \subset \R^n$. In view of condition $(C_1)$, this happens if and only if $E$ is of locally finite perimeter in the standard sense (corresponding to $\varphi(x,\xi)\equiv |\xi|$.)

In this paper we present existence, comparison, and uniqueness results for minimizers of the general least gradient problem 
\begin{equation}\label{LTVProb}
\inf_{v\in BV_f(\Omega)} \int_{\Omega} \varphi(x, Dv)
\end{equation}
where
$f\in C(\partial \Omega)$ and 
\[
BV_f(\Omega):=\{v\in BV(\Omega): 
\ \ \forall x\in \partial \Omega, \ \ 
\lim_{r\to 0} \ \esssup_{y\in \Omega, |x-y|<r} |f(x) - v(y)| = 0 \ \}.
%\ \ v \ \ \hbox{is continuous at}\ \ \partial \Omega \ \ \hbox{and},\ \ v=f \ \ \hbox{on} \ \ \partial \Omega\},
\]
%(We will sometimes say, slightly imprecisely, that ``$v$ is continuous at $\partial \Omega$, and $v=f$ on $\partial \Omega$",to describe the boundary condition in the definition of $BV_f(\Omega)$.)

\medskip

We will prove existence of a minimizer  in $BV_f(\Omega)$
of the general least gradient problem (\ref{LTVProb}),
as long as $\partial\Omega$ satisfies a positivity 
condition on a sort of generalized mean curvature related
to the integrand $\varphi$.
We refer to this as the {\em barrier condition}, and
we defer its statement until later (see Definition
\ref{def.barrier}, at the beginning of
Section \ref{S:exist}).
We will prove the following existence result. 
\begin{theorem}\label{existence.Theo}
Suppose that $\varphi: \R^n \times \R^n \rightarrow \R$ is a continuous function that satisfies $C_1-C_2$ in a bounded Lipschitz domain $\Omega\subset\R^n$.
If $\Omega$ satisfies the barrier condition  with respect to $\varphi$, then
for every $f\in C(\partial \Omega)$, the general least gradient problem (\ref{LTVProb}) has a minimizer in $BV_f(\Omega)$. 
\end{theorem}
In fact we prove something slightly stronger; see Remark \ref{slightlystronger}.

\medskip

For our comparison and uniqueness results we do not need to assume the barrier condition. On the other hand, we require stronger convexity and regularity conditions on $\varphi$
than we have so far assumed. In particular we will assume:\\

($C_3$) $\varphi\in W^{2,\infty}_{loc}$ away from $\{ \xi = 0\}$, and there exists $C>0$ such that $\varphi_{\xi_i \xi_j}(x,\xi) p^i p^j \ge C|p'|^2$
for all $\xi\in S^{n-1}$ and $p\in \R^n$, where $p' := p  - (p\cdot \xi)\xi$. \
\\

\medskip

($C_4$)
$\varphi$ and  $D_\xi \varphi$ are $W^{2,\infty}$ away from $\xi=0$, and 
there are positive constants $\rho$ and $\lambda$ such that

\begin{multline}
\varphi(x,\xi) + |D_\xi \varphi(x,\xi)| + |D^2_\xi \varphi(x,\xi)| + |D^3_\xi \varphi(x,\xi)| 
+
\rho|D_x D_\xi \varphi(x,\xi)|\\
+ \rho|D_x D^2_\xi \varphi(x,\xi)|
+ \rho^2|D_x^2 D_\xi \varphi(x,\xi)| \le \lambda \qquad\quad\mbox{ for all }x\in \Omega, \xi\in S^{n-1}.
\label{Fbounds}\end{multline}
These conditions are needed for a result,
due to Schoen and Simon \cite{SSA}, about partial regularity of
$\varphi$-minimizing sets, which we state as Theorem \ref{singularity-estim}.
This result will play a crucial role in our uniqueness proof in Section \ref{S:unique}. In addition, condition ($C_3$) is crucial also in Lemma \ref{HopfMax}.  
%One of the main ingredients in the proof of the above theorem is 

\medskip

The following theorem is our main uniqueness result. It also gives a stability estimate (with best constant) for the solutions with respect to errors in the boundary data.

\begin{theorem}\label{Comparison}
Let $\Omega \subset \R^n$ be a bounded Lipschitz domain with connected boundary, and assume $\varphi:\Omega \times \R^n\to \R$ satisfies $C_1-C_4$. Suppose  that 
$u_1, u_2$ minimize  \eqref{LTVProb} in $BV_{f_1}(\Omega)$ and  $BV_{f_2}(\Omega)$
respectively, for  $f_1, f_2\in C(\partial \Omega)$.
Then
\begin{equation}\label{uniqueness}
|u_2-u_1|\leq \sup_{\partial \Omega} |f_2-f_1| \ \ % \hbox{for}
 \ \ a.e. \ \  \mbox{in } \Omega.
\end{equation}
Moreover
\begin{equation}\label{comparison}
\mbox{ $u_2 \geq u_1$ a.e. in $\Omega$,  if $f_2 \geq f_1$ on $\partial \Omega$.}
\end{equation}
In particular, for every $f\in C(\partial \Omega)$, there is at most one minimizer of 
\eqref{LTVProb} in $BV_f(\Omega)$.
\end{theorem}

The following regularity result, valid only in low dimensions, is obtained by essentially the same arguments
as in the proof of Theorem \ref{Comparison}.
\begin{theorem} Let $\Omega \subset \R^n$ be a bounded Lipschitz domain with connected boundary, and assume $\varphi:\Omega \times \R^n\to \R$ satisfies $C_1-C_4$.
If $n=2$ or $3$ and $u$ minimizes  \eqref{LTVProb} in $BV_{f}(\Omega)$
for $f\in C(\partial \Omega)$, then $u$ is continuous in $\bar \Omega$.
\label{thm.cont}\end{theorem}

\medskip

For the special case $\varphi(x, \xi)=a(x)|\xi|$, the problem \eqref{LTVProb} is the weighted least gradient problem 

\begin{equation}\label{WLG}
\inf_{u\in BV_f(\Omega)} \int_{\Omega}a|Du|.
\end{equation}

In \cite{SZ}, %the authors studied minimization problem  (\ref{WLG}) with $a \equiv 1$. In particular, 
assuming a barrier condition on $\partial \Omega$ (see (3.1) and (3.2) in \cite{SZ}), the authors showed that  when $a\equiv 1$, problem (\ref{WLG}) admits a unique minimizer, and this minimizer is 
continuous in every dimension $n\ge 2$. Their proof is valid with very little change for problem \eqref{WLG} when $a$ is smooth enough and bounded below. It is not valid, however, for functionals $\varphi$
not arising from a Riemannian metric,
%satisfying hypotheses $C_1$ through $C_4$ above,
and it is not clear, and would be hard to determine, 
exactly how much regularity their proof requires of $a$. 
%Both of these limitations stem from the use made by \cite{SZ} of a maximum principle established in \cite{S2}.
The different approach we follow here yields uniqueness results and continuity
in low dimensions, under {\em sharp}
regularity hypotheses, for a larger class of problems, without any barrier condition on $\partial\Omega$.
We do not know whether, in the setting of Theorem \ref{Comparison} with $n\ge 4$, it is 
possible for a minimizer to be discontinuous.

We now briefly describe the conductivity imaging problem that leads to (\ref{WLG}). Let $\sigma(x)$ be a positive function that models inhomogeneous isotropic conductivity of a body $\Omega$.  If $u$ is the electric potential corresponding to the voltage $f$ on the boundary of $\Omega$, then $u$ solves the Drichlet problem
\[\nabla \cdot (\sigma \nabla u)=0, \ \ u|_{\partial \Omega}=f.\]
By Ohm's law, the corresponding current density is $J=-\sigma \nabla u$. Consider the inverse problem of determining $\sigma$ from knowledge of $|J|$ inside $\Omega$ (with a known $f$ prescribed on $\partial \Omega$).  Such internal data can be obtained using Magnetic Resonance Imaging \cite{joy89}.  It was first shown in \cite{NTT08} that the corresponding voltage potential $u$ is the unique solution of the weighted least gradient problem 
\begin{equation}
argmin \{ \int_{\Omega} a|D v|: \ \ u\in W^{1,1}(\Omega)\cap C(\bar{\Omega}) \ \ \hbox{and}\ \ u|_{\partial \Omega}=f\},
\end{equation}
with $a=|J|$ given. This uniqueness result has recently been extended to $u\in BV(\Omega)$ \cite{MNT13}. Once $u$ is determined (see \cite{MNTim} for a convergent numerical algorithm) the computation of $\sigma $ is straightforward. The uniqueness results in \cite{NTT08} and \cite{MNT13} assume that the weight $a$ is of the form $a=|J|$, as described above, but are valid for weights $a \in C^{\alpha}(\Omega)$ and allow $a$ to vanish in certain sets. The following direct consequence of Theorem \ref{Comparison} provides uniqueness, comparison, and stability results for general weights $a$ which are not necessarily of the form $a=|J|$. We do however need more restrictive assumptions on $a$. 

\begin{theorem}\label{uniqueness.Theo1}
Let $\Omega \subset \R^n$ be a bounded Lipschitz domain with connected boundary. Suppose $a\in C^{1,1}(\Omega)$ is positive and bounded away from zero, and $u_1, u_2\in BV(\Omega)$ 
minimize  (\ref{WLG}) in $BV_{f_1}(\Omega), BV_{f_2}(\Omega)$ respectively,
for some $f_1, f_2\in C(\partial \Omega)$.
Then
\begin{equation}\label{uniqueness2}
|u_2-u_1|\leq \sup_{\partial \Omega} |f_2-f_1| \ \ % \hbox{for}
 \ \ a.e. \ \  \mbox{in } \Omega.
\end{equation}
Moreover 
\begin{equation}\label{comparison2}
\mbox{ $u_2 \geq u_1$ a.e. in $\Omega$, if  $f_2 \geq f_1$ on $\partial \Omega$.}
\end{equation}
In particular,  (\ref{WLG}) has at most one minimizer in $BV_f(\Omega)$, and any
minimizer is continuous if $n\le 3$. 
\end{theorem}
In fact, in the setting of Theorem \ref{uniqueness.Theo1} (as well as Theorem \ref{uniqueness.Theo2} below)
minimizers are continuous for $\Omega \subset \R^n$ as long as $n \le 7$, see Remark \ref{X}.

\medskip

Our next result shows that the regularity assumption $a\in C^{1,1}(\Omega)$,
and hence also the regularity assumptions $C_3, C_4$  in Theorem \ref{Comparison}, 
are in a sense sharp. 

\begin{proposition} \label{non-unique-prop} For any $\alpha<1$, there exists a bounded smooth domain $\Omega \subset \R^n$  with connected boundary,  $f\in C(\partial \Omega)$, and a function $a \in C^{1,\alpha}(\Omega)$ with $\inf_{x\in \Omega}a(x)>0$ such the  weighted least gradient problem \eqref{WLG}
has infinitely many minimizers in $BV_f(\Omega)$.
\end{proposition}

In fact, since the weighted gradient functional is convex, 
 if uniqueness fails, then
there must be infinitely many minimizers.

The minimizers constructed in Proposition \ref{non-unique-prop} are all discontinuous, 
so in this sense the Proposition also shows that our regularity assumptions on $\varphi$ are  sharp
in Theorem \ref{thm.cont}.

It is also easy to see that all of our uniqueness and
comparison results
can fail if $\partial \Omega$ is not connected,
even if $a\in C^\infty(\bar \Omega)$. For example, if $\Omega$ is an annulus $B(2,0)\setminus \overline{B(1,0)}\subset \R^n$
and $a(x) = |x|^{1-n}$, then any function of the form $u(x) = g(|x|)$ minimizes \eqref{WLG}
with respect to its boundary data, as long as $g$ is monotone. Here and throughout the paper $B(r,x)$ denotes the open ball of radius $r$ centred at $x$.

\smallskip

Recently in \cite{HMN} authors presented a method for recovering the conformal factor of an anisotropic conductivity matrix in a known conformal class from one interior measurement.  Assume that the matrix valued conductivity $\sigma(x)$ is of the form
\[\sigma(x) =c(x)\sigma_0(x)\]
where $c(x)\in C^{\alpha}(\Omega)$ is a positive scalar valued function and $\sigma_0 \in C^{\alpha}(\Omega, Mat(n, \R^n))$ is  a known positive definite symmetric matrix valued function. In medical imaging  $\sigma_0$ can be determined using Diffusion Tensor Magnetic Resonance Imaging (see \cite{Ma} and the references therein). In \cite{MNT13} the authors showed that the corresponding voltage potential $u$ is the unique solution of the least gradient problem 
\[ argmin \{ \int_{\Omega} \varphi(x, Dv): \ \ u|_{BV(\Omega)}, \ \ u|_{\partial \Omega}=f \},\]
where $\varphi$ is given by 
\begin{equation}\label{varphi}
\varphi(x,\xi)=a(x)\left( \sum_{i,j=1}^{n}\sigma_0^{ij}(x)\xi _i \xi_j\right)^{1/2},
\end{equation}
\begin{equation}\label{aJJ}
a=\sqrt{\sigma_0^{-1} J \cdot J},
\end{equation}
and $J$ is the current density vector field generated by imposing the voltage $f$ at $\partial \Omega$.  Once $u$ is determined the function $c(x)$ can easily be calculated. 
This uniqueness result assumes that the weight $a$ is of the form (\ref{aJJ}) and it applies for weights $a \in C^{\alpha}(\Omega)$ that may vanish in certain sets (see \cite{HMN}). The following immediate corollary of Theorem \ref{Comparison} provides uniqueness, comparison, and stability results for general weights of the form (\ref{varphi}) with $a$ not necessarily of the form (\ref{aJJ}), but requires more restrictive assumptions on $a$.

\begin{theorem}\label{uniqueness.Theo2}
Let $\Omega \subset \R^n$ be a bounded Lipschitz domain with connected boundary, and assume $\varphi(x, \xi)$ is given by (\ref{varphi}), where $a\in C^{1,1}(\Omega)$ is positive and bounded away from zero  and $\sigma_0 \in C^{1,1}(\Omega, Mat(n,\R))$ satisfies 
\[ m |\xi|^2 \leq \sum_{i,j=1}^{n} \sigma_0^{ij} \xi_i \xi_j \leq M |\xi|^2, \ \ \hbox{for all} \ \ \xi \in \R^n,\]
for some $0<m, M<\infty $. 
If $u_1, u_2\in BV(\Omega)$ 
minimize  (\ref{WLG}) in $BV_{f_1}(\Omega), BV_{f_2}(\Omega)$ respectively, with
$f_1,f_2\in C(\partial \Omega)$,
then
\begin{equation}
|u_2-u_1|\leq \sup_{\partial \Omega} |f_2-f_1| \ \ % \hbox{for}
 \ \ a.e. \ \  \mbox{in } \Omega.
\end{equation}
Moreover
\begin{equation}
\mbox{ $u_2 \geq u_1$ a.e. in $\Omega$, if  $f_2 \geq f_1$ on $\partial \Omega$.}
\end{equation}

In particular,  for the class of $\varphi$ as described above, (\ref{LTVProb}) has at most one minimizer in $BV_f(\Omega)$,  and any
minimizer is continuous if $n\le 3$. 
\end{theorem}

The paper is organized as follows. In Section 2 we will present some preliminary results and basic facts about $\varphi$-total variation. Sections 3 and 4 are devoted to the proof of the existence and uniqueness results, respectively. In Section 5, we shall prove Proposition \ref{non-unique-prop} by constructing a one parameter family of minimizers of \eqref{WLG}. Finally, in Section 6, we provide a more convenient formulation of the barrier condition when the boundary of $\Omega$ is sufficiently smooth.

\section{Preliminary results}

In this section we develop some basic facts about $\varphi$-total variation, defined
in \eqref{varphi.def} above. These facts
are well-known for the usual notion of (isotropic, homogeneous) total 
variation, and we sketch some proofs here only to make it clear that the standard
arguments are still valid in the somewhat  more general setting that we consider here.  The paper \cite{AB} is a good reference for $\varphi-$total variation. 

It follows easily from \eqref{varphi.def1} that
\begin{equation}\
\int_A\varphi(x,D(u_1+u_2)) \le  \int_A\varphi(x,Du_1) + \int_A\varphi(x, Du_2)
\label{superadd}
\end{equation}
if $A$ is open, and hence if $A$ is any Borel set.\\

The definition \eqref{varphi.def} and the Fleming-Rishel coarea formula for $BV$ functions 
imply the following coarea formula for the $\varphi$-total variation.

\begin{proposition}[Remark 4.4 in \cite{AB}]
If $u\in BV_{loc}(\R^n)$ and $A\subset \R^n$ is a bounded Borel set, then
\begin{equation}\label{coarea}
\int_{A}\varphi(x, Du)= \int_{-\infty}^{+\infty} P_{\varphi}(X_t; A)dt,
\end{equation}
where $X_t:=\{x \in \Omega : u(x)>t\}$ and $P_\varphi$ is as defined in equation (\ref{perimeter}).  
\end{proposition}
It is a straightforward consequence of the coarea 
formula that for any $u\in BV(\Omega)$ and $\lambda\in \R$, 
if we write $u_1 :=  \max(u-\lambda, 0)$ and $u_2 = u - u_1$, then
\begin{equation}
\int_{\Omega} \phi(x, Du) = \int_\Omega \phi(x, D u_1) + \int_\Omega \phi(x, Du_2).
\label{additive}\end{equation}

\begin{lemma}
Let $A\subset \R^n$ be a Borel set and $E_1, E_2\subset \R^n$ be of locally finite perimeter with respect $\varphi$. Then
\begin{equation}\label{areMinIneq}
P_{\varphi}(E_1 \cup E_2; A)+P_{\varphi}(E_1 \cap E_2; A) \leq P_{\varphi}(E_1; A)+P_{\varphi}(E_2; A). 
\end{equation}
\end{lemma}

{\bf Proof:}
We apply \eqref{additive} with $u= \chi_{E_1}+\chi_{E_2}$ and $\lambda = 1$.
Then $u_1 = \chi_{E_1\cap E_2}$,  and $u_2 = 
\chi_{E_1\cup E_2}$.
It follows that 
\begin{align*}
\int_{A} \varphi (x, D\chi_{E_1 \cup E_2})ds+\int_{A} \varphi(x,D\chi_{E_1 \cap E_2})ds
&
\overset{\eqref{additive}}=
\int_{A} \varphi(x, D(\chi_{E_1}+\chi_{E_2}))\\
&
\overset{\eqref{superadd}}\le
\int_{A} \varphi(x, D\chi_{E_1})
+
\int_{A} \varphi(x, D\chi_{E_2}).
\end{align*}
Rewriting in terms of $P_\varphi$ yields \eqref{areMinIneq}.
\hfill $\Box$
\\

\begin{definition}
(i) We say that a function $u\in BV(\R^n)$ is $\varphi$-total variation minimizing in a
set $\Omega\subset \R^n$ if 
\[
\int_{R^n} \varphi(x, Du) \le 
\int_{R^n} \varphi(x, Dv) \qquad\mbox{ for all $v\in BV(\R^n)$ such that $u=v$ a.e. in $\Omega^c$}.
\]
(ii) Similarly, we say that $E\subset \R^n$ of finite perimeter is $\varphi$-area minimizing in $\Omega$
if
\[
P_{\varphi}(E) \leq P_{\varphi}(F) 
\quad\mbox{for all $F\subset \R^n$ such that  $F \cap \Omega^c=E\cap \Omega^c$\ a.e.. }
\]
\end{definition} 
We emphasize that in the definitions above, $u-v$ is not required to have compact support in $\Omega$, and $E\setminus F$ is not required to be
compactly contained in 
$\Omega$.
\\

If $v\in BV(\R^n)$ and $\Omega$ is an open set with Lipschitz boundary, we will write
$v^+$ and $v^-$ to denote the outer and inner trace of $v$ on $\partial \Omega$.
Recall that these are functions in $L^1(\partial \Omega; \mathcal H^{n-1})$,
characterized by the fact that for $\mathcal H^{n-1}$ almost every  $x\in \partial \Omega$,
\begin{equation}\label{traces}
\lim_{\rho\to 0} \frac 1{\rho^n}\int_{B_\rho(x)\setminus \Omega}|v^+(x) - v(y)| \ dy \  = \  
\lim_{\rho\to 0} \frac 1{\rho^n}\int_{B_\rho(x)\cap \Omega}|v^-(x) - v(y)| \ dy = 0.
\end{equation}

\begin{lemma}Let $\Omega\subset \R^n$ be bounded and open, with Lipschitz boundary.
Given $g\in L^1(\partial \Omega ; \mathcal H^{n-1})$, define 
\[
I_\varphi(v ; \Omega, g) := \ \int_{\partial \Omega} \varphi(x, \nu_{\Omega}) |g - v^-|
d\mathcal H^{n-1}
+
\int_{\Omega}\varphi(x,Dv).
\]
where $\nu_\Omega$ denotes the outer unit normal to $\Omega$.
Then $u\in BV(\R^n)$ is $\varphi$-total variation minimizing in $\Omega$ 
if and only if $u|_\Omega$ minimizes
$I_\varphi ( \, \cdot \, ; \Omega, g)$ for some $g$, and moreover $g = u^+$.
\label{lem:reformulate}\end{lemma}

{\bf Proof:}
We recall some basic properties of traces.
First, if $v\in BV(\R^n)$ then $v^+$ and $v^-$ are $\mathcal H^{n-1}$ integrable
on $\partial \Omega$, and conversely, for every $g\in L^1(\partial \Omega;\mathcal H^{n-1})$
there exists some $v\in BV(\R^n)$ such that $g = v^+$ say.
Second, we note that
\begin{equation}
\int_{\partial \Omega} \varphi(x, Dv) = \int_{\partial \Omega} \varphi(x, v^v) |Dv|
=
\int_{\partial \Omega} \varphi(x, \nu_{\Omega} )|v^+-v^-| d\mathcal H^{n-1}.
\label{boundary.meas}\end{equation}
To see this, note that $|Dv|$ can only concentrate on a
set of dimension $n-1$ if that set is a subset of the jump set of $v$,
so \eqref{boundary.meas} follows from standard descriptions of
the jump part of $Dv$.

Now if $u, v\in BV(\R^n)$ satisfy $u=v$ a.e. in $\Omega^c$, then it follows from \eqref{varphi.def} that $\int_{\bar \Omega^c} \varphi(x, Du) = \int_{\bar \Omega^c} \varphi(x, Dv)$. In
addition, $u^+ = v^+$, so
using  \eqref{boundary.meas}
we deduce that
\[
\int_{\R^n}\varphi(x, Du) - \int_{\R^n} \varphi (x,Dv)
 \ = \ I_\varphi(u; \Omega, u^+) - I_\varphi(v ; \Omega, u^+).
\]
The lemma easily follows.
\hfill $\Box$
\\

\begin{lemma}\label{areaMinCor}
Let $\Omega$ be a bounded Lipschitz domain and let $E_1, E_2\subset \R^n$ be area minimizing in $ \Omega$.  If $E_1 \cap  \Omega^c  \subset E_2 \cap  \Omega^c$, % and $E_1\setminus E_2 \subset \subset \Omega$, 
then $E_1\cap E_2$ and $E_1 \cup E_2$ are area minimizing in $ \Omega$. 
\end{lemma}

{\bf Proof:} For $i=1,2$, let $\mathcal A_i := \{ F \subset \R^n : F\cap  \Omega^c = E_i \cap  \Omega^c \}$.
Our hypotheses imply that $E_1\cap E_2 \in \mathcal A_1$ and $E_1\cup E_2\in \mathcal A_2$,
so it suffices to show that 
\[
P_\varphi(E_1\cap E_2) \le \inf_{F\in \mathcal A_1} P_\varphi(F) = P_\varphi(E_1),
\ \ \ 
P_\varphi(E_1\cup E_2) \le \inf_{F\in \mathcal A_2} P_\varphi(F) = P_\varphi(E_2).
\]
Since the opposite inequalities hold, these follow directly from \eqref{areMinIneq}.
\hfill $\Box$
\\

\begin{theorem}\label{areMinTheo}
Let $\Omega \subset \R^n$ be a bounded Lipschitz domain and $u\in BV(\R^n)$ be $\varphi$-total variation minimizing in $\Omega$.
Let 
\[E_{\lambda}=\{x\in \R^n: \ \ u(x)\geq \lambda \}. 
\]
Then $E_\lambda$ is $\varphi$-area minimizing in $\Omega$ for  every $\lambda$. 
\end{theorem}

A similar result used in earlier work
about the case $\varphi(x,\xi) = |\xi|$ (see for example \cite{SZ}) uses a somewhat different notion of minimizing, in which
only perturbations with compact support in $\Omega$ are allowed. For the proof we will use the following lemma.

\begin{lemma}\label{prop1}
Assume $u_k$ is $\varphi$-total variation minimizing in $\Omega$
for $k \geq 1$ and
\[
\mbox{ $u_k\to u$ in $L^1(\Omega)$,  and  $u^\pm_k\to u^\pm$
in $L^1(\partial \Omega; \mathcal H^{n-1})$. }
\]
Then $u$ is $\varphi$-total variation minimizing in $\Omega$.
\end{lemma}

{\bf Proof:}
It follows  from \eqref{varphi.def1} via quite standard arguments
that 
%Let $g := u^+$ on $\partial \Omega$.
%It suffices to show that for any $v\in BV(\Omega)$,
%\[
%I_\varphi(u;\Omega,g) \le I_\varphi(v ; \Omega,g).
%\]
\begin{equation}
\int_{\Omega}\varphi(x, Du) \le 
\liminf_k
\int_{\Omega}\varphi(x, Du_k),
\label{wlsc}\end{equation}
and this, with  the $L^1$ convergence of the traces, implies that
\begin{equation}\label{limit1}
I_\varphi(u; \Omega, u^+) \le \liminf_{k\to\infty} I_\varphi(u_k; \Omega, u_k^+).
\end{equation}
Now for any $v\in BV(\R^n)$ such that $u=v$ a.e. in $\Omega^c$, 
\begin{align*}
I_\varphi(u_k; \Omega, u_k^+) 
&\le
I_\varphi(v; \Omega, u_k^+)\\
&\le
I_\varphi(v; \Omega, u^+) +\int_{\partial \Omega}\varphi(x,  \nu_\Omega)
|u^+ - u_k^+| \ d\mathcal H^{n-1}\\
&\le
I_\varphi(v; \Omega, u^+) +\alpha^{-1} \int_{\partial \Omega}
|u^+ - u_k^+|\ d\mathcal H^{n-1}
\end{align*}
using Lemma \ref{lem:reformulate}, the minimality of $u_k$ and
standing assumption $C_1$. It follows from this, \eqref{limit1},
and again the $L^1$ convergence of the traces that
$I_\varphi(u;\Omega, u^+) \le 
I_\varphi(v;\Omega, u^+) $,
which proves the proposition.
\hfill $\Box$

%\begin{remark}
%The hypothesis that $u_k^+\to u^+$ in $L^1(\partial \Omega)$
%is necessary for the above proposition, as is shown by the example (in 2 dimensions):
%\[
%5u_k = \chi_{\{ x\in \R^2 : |x| > 1+\frac 1k\}},
%\qquad
%u = \chi_{\{ x\in \R^2 : |x| > 1\}},
%\]
%and $\Omega$ the unit disk. For every $k$, $u_k$ is minimizing in $\Omega$
%with respect to the usual total variation
%(%zero boundary data, and zero in $\Omega$) but the limiting function
%(boundary data 1 but zero in $\Omega$)
%is certainly not a minimizer. I don't know whether we want to mention this. I also don't know whether the hypothesis $u^-_k\to u^-$
%is needed --- I can't think of any counterexample, and it kind of seems like it should be superfluous, but I couldn't see a way to get rid of it either.
%\end{remark}

{\bf Proof of Theorem \ref{areMinTheo}:} Our argument is modelled on the proof of Theorem 1 in %Bombieri , de Giorgi, Giusti 
\cite{BdGG}.
 
For $\lambda\in \R$, let $u_1 =  \max(u-\lambda, 0)$, $u_2 = u - u_1$.
Let $g\in BV(\R^n)$ with $supp(g)\subset \bar{\Omega}$.
Since $u$ is a minimizer,
\begin{eqnarray*}
\int_{\Omega}\varphi(x,Du_1)+\int_{\Omega} \varphi(x, Du_2) 
&
\overset{\eqref{additive}} =
&\int_{\Omega}\varphi (x, Du)\\
&\leq & \int_{\Omega}\varphi(x,D(u+g))\\
&\overset{\eqref{superadd}} \leq & \int_{\Omega} \varphi(x,D(u_1+g))+\int_{\Omega}\varphi(x,Du_2). 
\end{eqnarray*}
Hence $u_1$ is a minimizer. Repeating the same argument, one verifies that
\[
\chi_{\epsilon, \lambda}:= \min(1, \frac 1 \epsilon u_1  )
= \begin{cases}
0&\mbox{ if }u\le \lambda\\
\epsilon^{-1}(u-\lambda)&\mbox{ if }\lambda \le u \le \lambda+\epsilon\\
1&\mbox{ if }u\ge \lambda+\epsilon
\end{cases}
\]
is also a minimizer of (\ref{varphi.def}).  
It is clear that for a.e. $\lambda\in \R$,
\begin{equation}\label{goodlambda}
\mathcal L^n( \{x\in \Omega : u(x) = \lambda\})
=
\mathcal H^{n-1}( \{x\in \partial \Omega : u^\pm(x) = \lambda\})
= 0,
\end{equation}
and it is straightforward to check, using \eqref{traces}, 
that if \eqref{goodlambda} holds, then
\begin{equation}
\chi_{\epsilon,\lambda}\to \chi_\lambda := \chi_{E_\lambda}
\mbox{ in }L^1_{loc}(\R^n),\qquad
\chi^{\pm}_{\epsilon,\lambda}\to \chi_\lambda^\pm
\mbox{ in }L^1(\partial \Omega;\mathcal H^{n-1}).
\label{approx.step}\end{equation}
Thus Lemma \ref{prop1} implies that 
$\chi_{E_\lambda}$ is $\varphi$-total variation minimizing
in $\Omega$, and hence that
$E_\lambda$ is 
$\varphi$-area minimizing in $\Omega$.

If $\lambda$ does not satisfy \eqref{goodlambda}, then let $\lambda_k$ be an increasing sequence
such that $\lambda_k \to \lambda$ and $\lambda_k$ satisfies \eqref{goodlambda} for every $k$.
Then one can check that 
\[
\chi_{\lambda_k}\to \chi_\lambda \quad \mbox{ in }L^1_{loc}(\R^n),
\quad
\chi^{\pm}_{\lambda_k}\to \chi_\lambda^\pm
\mbox{ in }L^1(\partial \Omega;\mathcal H^{n-1}).
\]
as $k\to \infty$, so it again follows from Proposition \ref{prop1}
that $E_\lambda$ is 
$\varphi$-area minimizing in $\Omega$.
\hfill $\Box$

\begin{definition}
Let $E \subset \R^n$. A point $x\in \partial E$ is called a regular point if there exists $\rho>0$ such that $\partial E \cap B(x,\rho)$ is a $C^2$ hypersurface.  We denote the set of all regular points of $\partial E$ by $reg (\partial E)$. We will say that $x$ is a singular point of $\partial E$ if $x \in$sing$(\partial E)$, where
\[sing(\partial E)= \partial E \setminus reg (\partial E).\]
\label{Reg.def}\end{definition}
If $E$ is a measurable subset of $\R^n$, we will write 
\begin{equation}
E^{(1)}:= \{x\in \R^n :
\lim_{r\to 0} \frac {\calH^n(B(r,x)\cap E)}{\calH^n(B(r))} = 1 \}.
\label{E0E1}\end{equation}

The following singularity estimate is due to Schoen, Simon, and Almgren \cite[Theorem I.3.1 and Corollary I.3.2]{SSA},  and plays a crucial role in our uniqueness proof in Section \ref{S:unique}.

\begin{theorem}\label{singularity-estim}
Let $\Omega \subset \R^n$.  Suppose $\varphi: \R^n \times \R^n\rightarrow \R$ satisfies $C_1-C_4$.
%, and $\frac{\partial \varphi}{\partial \xi}$ exists  and is $C^{1,1}$ away from $\xi=0$. 
If $E$  is $\varphi$-area minimizing in $\Omega$, 
%and $E$ is normalized in the sense \eqref{normalized.def}, 
then
\begin{equation}\label{DimSingSet}
\left\{ \begin{array}{ll}
\mathcal{H}^{n-3}(\hbox{sing}(\partial E^{(1)}) \cap \Omega)< \infty \ \ &\hbox{if} \ \ n\geq 4\\
\hbox{sing}(\partial E^{(1)}) \cap \Omega =\emptyset, \ \ &\hbox{if} \ \ n \le 3.
\end{array} \right.
\end{equation}
\end{theorem}

One reason for considering $\partial E^{(1)}$ rather than $\partial E$ is that the 
former is insensitive to modifications of $E$ on sets of measure $0$.

\begin{remark}\label{rem:closed}
In \cite{SSA}, the set that we have identified as $\partial E^{(1)}$, for a $\varphi$-area minimizing set  $E$, is described in a slightly different way,
but standard facts about BV functions imply that our description
is equivalent in the context of the theorem.
\end{remark}

\begin{remark}\label{rem:SSA}
Some ambiguous wording in \cite{SSA} suggests that  Theorem \ref{singularity-estim} might
require that $\varphi \in W^{3,\infty}_{loc}$ away from $\{\xi = 0\}$, in addition
to the hypotheses in  $C_4$. However,  inspection of
the proof shows that $C_4$ is all the regularity that is needed for $\varphi$.
Indeed, regularity of $\varphi$ is used in the proof of %Theorem \ref{singularity-estim} in 
\cite[Theorem I.3.1]{SSA} in the following ways:
\begin{enumerate}

\item to ensure $C^{2,\gamma}_{loc}$ regularity for weak
solutions  $w: B^{n-1}_r\to \R$ of equations of the form
\[
- \sum_{i=1}^{n-1} \partial_{x_i} (\varphi_{\xi_i}(x,w ,-Dw,1)) - \varphi_{x_n}(x,w,-Dw,1) = 0 , \qquad
x\in B^{n-1}_r.
\]

\item as a hypothesis for basic $\epsilon$-regularity results for
$\varphi$-area minimizing currents.

\item for maximum principle arguments such as those in Lemma \ref{HopfMax} below

\item in estimates such as those in the proof of \cite[Lemma I.2.5]{SSA}.

\end{enumerate}
For the first point listed above, it is rather standard that the regularity assumed in $C_4$ suffices, together with the structural conditions $C_1$ - $C_3$, and for the second, it is proved in the
reference \cite{SS} cited in \cite{SSA} that these hypotheses are enough.
For the maximum principle argument, it is clear from
our proof of Lemma \ref{HopfMax} below (which is just a slightly more detailed
version of an argument from \cite{SSA}) that $C_1 - C_3$ suffice. As well, all the estimates
in \cite{SSA} involve only the derivatives appearing in $C_4$, and yield constants
that depend only on the dimension $n$ and the constants $\rho, \lambda$ in $C_4$. 
\end{remark}

\section{Existence}\label{S:exist}

In this section we prove Theorem \ref{existence.Theo}. First we give a precise statement of the
geometric condition that is a main hypothesis of the theorem.

\begin{definition} Let $\Omega \subset \R^n$ be a bounded Lipschitz domain and $\varphi: \Omega \times \R^n \rightarrow \R$ is a continuous function that satisfies $C_1-C_2$. We say that $ \Omega$ satisfies the  barrier condition if for every $x_0\in \partial \Omega$ and $\epsilon>0$ sufficiently small, if $V$
minimizes $P_{\varphi}(\,\cdot\,; \R^n)$ in
\begin{equation}\label{SBC-Set}
\{W \subset \Omega : W \setminus  B(\epsilon, x_0)=\Omega \setminus B(\epsilon, x_0)\},
\end{equation}
then
\[\partial V^{(1)} \cap \partial \Omega \cap B(\epsilon, x_0)= \emptyset.\]
\label{def.barrier}\end{definition}

In the case $\varphi(x,\xi)= |\xi|$, the barrier condition is equivalent, at
least for smooth sets, to the one introduced in \cite{SZ}.

%{The barrier condition is necessary, as well as sufficient for existence of a solution to the minimization problem \eqref{LTVProb}. Indeed, if the barrier condition fails at some point $x_0\in \partial \Omega$, then \eqref{LTVProb} has no solution for boundary  data $f(x) = \mbox{dist}(x_0, x)$. }

A convenient
interpretation of the barrier condition, if  $\partial \Omega$ is sufficiently smooth,
is provided by the following result. See also Remark \ref{R.below}  below.
\begin{lemma}
Assume that $\partial \Omega$ is $C^{2}$
and that $\varphi$ satisfies $C_1-C_3$. Define the signed distance $d(\cdot)$ to
$\partial \Omega$ by
\[
d(x) := \begin{cases}
\operatorname{dist}(x, \partial \Omega) &\mbox{ if }x\in \Omega\\
-\operatorname{dist}(x, \partial \Omega) &\mbox{ if not}.
\end{cases}
\]
Then $\Omega$ satisfies the barrier condition if
\begin{equation}
-\sum_{i=1}^n \partial_{x_i} \varphi_{\xi_i} (x, Dd(x)) >0 \quad\mbox{ on a dense subset of 
$\partial \Omega$}.
\label{barrier.suff}\end{equation}
\label{L.barrier}\end{lemma}

The proof of Lemma \ref{L.barrier}
is given in  Section \ref{S:barrier}.
Although we do not prove it, if $\varphi$ and $\Omega$
satisfy the above hypotheses, then the barrier 
condition is in fact {\em equivalent} to 
\eqref{barrier.suff}.

\begin{remark}\label{R.below}
If $\partial \Omega$ is defined (locally) as a graph $x_n=w(x')$ of a $C^2$ function $w$ (so that locally $\Omega=\{(x',x_n) \in \R^{n-1}\times \R :  x_n>w(x')\}$) then an equivalent formulation of (\ref{barrier.suff}) can be seen from the following equality 
\begin{equation}\label{equ}
\sum_{i=1}^n \partial_{x_i}\varphi_{\xi_i}(x, Dd)
=
\sum_{j=1}^{n-1} \partial_{x_j}\varphi_{\xi_j}(x', w, -Dw, 1)
+\varphi_{x_n}(x', w, -Dw, 1),
\end{equation}
at points $x=(x', w(x'))$ on $\partial \Omega$. For the convenience of the readers we include a proof of (\ref{equ}) in Section 6. 
\end{remark}

\medskip

Our main use of the barrier condition is the following technical lemma.

\begin{lemma}\label{SBC-Lemma}
Let $\Omega \subset \R^n$ be a bounded Lipschitz domain that satisfies the  barrier condition 
with respect to $\varphi$, and assume that $E \subset \R^n$ is $\varphi-$area minimizing in $\Omega$.
Then
 \[
 \{x\in \partial \Omega \cap \partial E^{(1)}: \ \ B(\epsilon, x) \cap \partial E^{(1)}\subset \bar{\Omega} \ \ \hbox{for some} \ \ \epsilon>0\}=\emptyset.
 \] 
\end{lemma}
{\bf Proof:} Assume there exists $x_0\in \partial \Omega \cap \partial E^{(1)}$ such that  $B(\epsilon, x_0) \cap \partial E^{(1)}\subset \bar{\Omega}$ for some $\epsilon>0$ and let $V$ be a minimizer of $P_{\varphi}(\,\cdot\,; \R^n)$ in (\ref{SBC-Set}).  Existence of such a set $V$ is standard,
for reasons discussed in the proof of Theorem  \ref{existence.Theo} below.

Then it follows from Lemma \ref{areaMinCor} that $V'=V\cup (E \cap \Omega)$ also minimizes $P_{\varphi}(\,\cdot\,; \R^n)$ in (\ref{SBC-Set}).
It is easy to see that  $x_0\in \partial {V'}^{(1)}$. 
Hence 
\[
x_0\in \partial {V' }^{(1)}\cap \partial \Omega \cap B(\epsilon, x_0)\neq  \emptyset. 
\]
This contradicts the  barrier condition and finishes the proof. 
\hfill $\Box$

\medskip

Since $f \in C(\partial \Omega)$, it can be extended to a function in  $C(\Omega^c)$ and throughout the paper we will denote this extension to $C(\Omega^c)$ by $f$ again. Now we are ready to prove our existence result. \\

{\bf Proof of Theorem \ref{existence.Theo}:}  
Without loss of generality we may assume that $f\in BV(\R^n)$,
since basic trace theorems guarantee that every $\calH^{n-1}$
integrable function on $\Omega$ is the trace of some
(continuous) function in $BV(\Omega^c)$.
Define
\[
\calA_f := \{v\in BV(\R^n): \ \ v=f \ \ \hbox{on}\ \ \Omega^c \},
\]
and note that $BV_f(\Omega) \hookrightarrow \calA_f$, 
in the sense that any element $v$ of $BV_f(\Omega)$
is the restriction to $\Omega$ of a unique element of 
$\calA_f$. Thus it suffices to prove that 
the functional
\[
F(v):=\int_{\R^n}\varphi(x,Dv) dx.
\]
has a minimizer $u\in \calA_f$, and that $u$ 
can be identified with an element of $BV_f(\Omega)$.

Existence of a minimizer $u$ is standard, as
$F$ is coercive in $BV(\R^n)$ (a consequence of $C_1$) and
weakly lower semicontinuous, as already noted in  \eqref{wlsc},
and because $BV(\R^n)\hookrightarrow L^1_{loc}$.

We next use the barrier condition to show  that $u\in BV_f(\Omega)$.
If not, there exists some $x\in \partial \Omega$
and $\delta>0$ such that
\begin{equation}
 \esssup_{y\in \Omega, |x-y|<r}\big( f(x) - u(y))\ge \delta\qquad
\mbox{ or }  \ \ \esssup_{y\in \Omega, |x-y|<r}\big(u(y) - f(x)) \ge \delta
\label{ess.alt}\end{equation}
for every $r>0$.
Assume that the latter condition holds. 
It follows
from this, along with the definition of $\calA_f$ and the continuity of $f$, 
that $x\in  \partial E^{(1)}$ for $E := E_{ f(x) + \delta/2}$. 
Recall that  $E$ is $\varphi$-area minimizing in $\Omega$
by Theorem \ref{areMinTheo}.
However, since $f$ is continuous in $\Omega^c$ and
$u\in \calA_f$, it is clear that
$u< f(x) + \delta/2$ in $B(\e, x)\setminus \Omega$ for all sufficiently
small $\e$.
But  Lemma \ref{SBC-Lemma}
shows that this is impossible.

If the first alternative in \eqref{ess.alt} holds, then we
define $E :=\{ y\in \R^n : u(y) \le f(x) - \delta/2\}$
and find in the same way that $u\in BV_f(\Omega)$.

Finally, note that if
$v\in BV_f$, then  the inner and outer traces
of $v$ both equal $f$ at every point of $\partial \Omega$,
and so it follows
from \eqref{boundary.meas} that
that $|Dv|(\partial \Omega) = 0$. 
Hence, if $v\in BV_f(\Omega)$, then
\begin{align*}
F(v)&=
\int_\Omega \varphi(x, Dv) +
\int_{\partial \Omega} \varphi(x, Dv) +
\int_{\R^n\setminus \bar \Omega} \varphi(x, Dv)\\
&=
\int_\Omega \varphi(x, Dv) +
\int_{\R^n\setminus \bar \Omega} \varphi(x, Df)
\end{align*}
So the fact that $u\in BV_f(\Omega)$ minimizes $F$ in $\calA_f$ 
implies that it minimizes  \eqref{LTVProb} in $BV_f(\Omega)$.

\hfill $\Box$

\begin{remark}\label{slightlystronger}
The above arguments show that if $u\in BV_f(\Omega)\hookrightarrow \calA_f$
minimizes \eqref{LTVProb} in $BV_f(\Omega) $, then
it also is $\varphi$-total variation minimizing in $\Omega$, so that
$F(u) \le F(v)$ for all $v\in \calA_f\supsetneq BV_f(\Omega)$. This conclusion does not
require the barrier condition.
\end{remark}

% If $u$ is not continuous at $\partial \Omega$, then (\ref{continuity}) does not hold for some $\lambda \in \R$. Since both $u$ and $-u$ minimize $F$ relative to their own boundary data, without loss of generality we may assume that  \[u(\partial E_\lambda \cap \partial \Omega) \not \subset [\lambda, \infty). \] Thus there exists $x_0\in \partial \Omega \cap \partial E_\lambda$ with $u(x_0)<\lambda$. Since $u$ is continuous on $\Omega^c$, $B(\epsilon, x_0)\cap \partial E_\lambda \subset \bar{\Omega}$ for some $\epsilon>0$.  In view of Lemma \ref{SBC-Lemma} this is a contradiction. Hence $u$ is continuous at $\partial \Omega$ and the proof is complete.

\section{Uniqueness and continuity}\label{S:unique}

In this section we prove the comparison principle and uniqueness result stated in Theorem \ref{Comparison}. We shall use the following results from dimension theory.  The proofs can be found in \cite{HW} (Chapter IV).

\begin{proposition}\label{openProp}
Let $U$ be an open set in a connected $k$-dimensional manifold  which is neither empty nor dense. Then the topological dimension of $\partial U$ is $k-1$.
\end{proposition}

\begin{proposition}\label{disconnect-prop}
A connected $k$-dimensional manifold can not be disconnected by a subset of dimension $k-2$.
\end{proposition}

For a proof of the following proposition see \cite{HW} (Chapter VIII, $\S 4$ ).
\begin{proposition}\label{HTdim}
Let $X$ be a metric space. Then the Hausdorff dimension of $X$ is bounded below by its topological dimension. 
\end{proposition}

Let $E \subset \R^n$  be a  $\varphi$-area minimizing set and $y_0$ be a regular point of $\partial E^{(1)}$. Then  for $\rho$ sufficiently small,
we can arrange, after a suitable choice of coordinates, that
\[\partial E^{(1)} \cap B(\rho, y_0)=\{(y,w(y)): \ \ y\in A\},\]
for some open $A\subset \R^{n-1}$ and $w \in C^2(A)$.  Moreover,  by rewriting
$\int \varphi (x,D \chi_E)$ in terms of $w$ and computing the first variation,
we find that $w$ satisfies 
\begin{equation}\label{HopfEq}
\varphi_{x_n}(y,w, -Dw, 1)+\sum_{i=1}^{n-1} \frac{d}{dx_i} (\varphi_{\xi_i}(y, w, -Dw, 1))=0, \ \ y\in A.
\end{equation}

The following strong maximum principle is a standard
consequence of basic elliptic regularity results.

\begin{lemma} \label{HopfMax} Suppose $\varphi$ satisfies $C_1-C_3$.
Assume also that $w_1$ and $w_2$ are  $C^2$ solutions of (\ref{HopfEq}) on a 
$(n-1)$-dimensional
ball $B(y_0,\rho)$ such that  $w_1 \le w_2$, and that $w_1 = w_2$ at some point
in $B(y_0, \rho)$. Then $w_1=w_2$ on $B(y_0,\rho)$. 
\end{lemma}

Note that $C_4$ is not needed here. The failure of Lemma \ref{HopfMax} for less smooth
integrands $\varphi$ is the mechanism behind the counterexamples presented in the next section. 
In order to make it clear where the regularity assumptions on $\varphi$ are used
we therefore present some details of the proof, which however is quite standard (see {\em e.g. }\cite[Lemma 2.4]{SSA}).

\medskip

{\bf Proof:} Let $w=w_2-w_1$. If we write $\mathcal D(w)$ for the left-hand side of \eqref{HopfEq},
then by rewriting the identity $0 =\mathcal D(w_2) - \mathcal D(w_1) =  \int_0^1 \frac d{ds}\mathcal D(s w_2 + (1-s)w_1)\ ds$,
we find that $w$ satisfies the equation
\begin{equation}
\sum_{i,j=1}^{n-1} \frac{\partial}{\partial x_i} \left(-\langle \varphi_{\xi_i \xi_j}\rangle w_{x_j} + \langle \varphi_{\xi_i x_n}\rangle w \right)
+\langle\varphi_{x_n \xi_i} \rangle w_{x_i} 
- \langle\varphi_{x_n x_n}\rangle w \ = \ 0 , 
\label{linear.diff}\end{equation}
where we use the notation
\[
\langle g \rangle (y) := \int_0^1 g(y, w^s(y), -Dw^s(y), 1) \ ds, \qquad w^s = s w_2 + (1-s)w_1.
\]
It follows from our assumptions that all the coefficients in the above equation are 
$L^\infty$, and that $\langle \varphi_{\xi_i \xi_j}\rangle$ is positive definite.
Since $w\ge 0$ and $w$ vanishes at some point in
$B(y_0,\rho)$, it follows from Moser's Harnack inequality, which is valid for
equation \eqref{linear.diff}, that $w \equiv 0$ in $B(y_0,\rho)$. $\Box$\\

We shall also  need the following lemma. 

\begin{lemma}\label{intersectLem} Let $\Omega$ be a bounded Lipschitz domain with connected boundary and assume that $E \subset \R^n$ is $\varphi-$area minimizing in $\Omega$.  If $R$ is a nonempty connected component of $reg(\partial E^{(1)})\cap \Omega$, then
$\bar R \cap \partial \Omega \neq  \emptyset$.
\end{lemma}

Recall that $reg(\partial E^{(1)})$ denotes the regular part of $\partial E^{(1)}$ (see Definition \ref{Reg.def}).
The idea of the proof is that if the conclusion fails, then we could modify $E$ in a way that
decreases its $\varphi$-perimeter without changing $E\cap \Omega^c$,
either by deleting a component  of $E$ or by ``filling in a hole"; this would contradict the 
minimality of $E$. We defer the full proof to the end of the section.

\smallskip

%We shall need the following results from dimension theory.  The proofs can be found in \cite{HW} (Chapter IV).

We will deduce Theorems \ref{Comparison} and  \ref{thm.cont}
from the following geometric comparison principle,
which is of interest in its own right.

\begin{theorem}[Comparison principle for $\varphi$-area minimizing sets]\label{MaxPrinc}
Let $\Omega \subset \R^n$ be a bounded Lipschitz domain with connected boundary, and suppose  that $\varphi$ satisfies $C_1-C_4$. Assume that 
$E_1, E_2\subset \R^n$ are $\varphi$-area minimizing sets 
in $\Omega$.
and also  that 
\begin{equation}\label{CC}
E_1 \setminus \Omega \subset \subset E_2 \setminus \Omega.
\end{equation}

%\[
%\hbox{for}\ \ i=1, 2,\quad
% \mbox{$E_i \subset \R^n $ minimizes $P_{\varphi}(\,\cdot\,; \R^n)$ in}\ \
%\{W \subset \R^n: W\backslash \Omega = E_i \backslash \Omega\}
% \]
%and also  that 
If  $\Omega$ satisfies the barrier condition, or if
%\begin{equation}\label{EmptyIntersec}
%\partial E_1^{(1)}  \cap \partial \Omega  \subset E_2^{(1)}, \qquad\qquad
%\partial E_2^{(1)} \cap \partial \Omega \subset \partial \Omega \setminus  \overline {E_1^{(1)}}
%\end{equation}
\begin{equation}\label{EmptyIntersec}
\partial E_1^{(1)}  \setminus E_2^{(1)}\subset \Omega, 
 \quad\qquad\mbox{ and }\qquad
\partial E_2^{(1)} \cap  \overline E_1^{(1)} \subset \Omega,
\end{equation}
then
\[
E_1^{(1)} \subset E_2^{(1)}.
%\calH^n(E_1 \setminus E_2)=0.
\]
Moreover,  if $n\le 3$ then $E_1^{(1)}\subset \subset E_2^{(1)}$.
\end{theorem}

\begin{remark}
When applying the above theorem in the proof of our main uniqueness result, Theorem \ref{Comparison}, assumption (\ref{EmptyIntersec}) will be satisfied as a consequence of the sense in which the minimizers $u_1, u_2$ assume their boundary values $f_1, f_2$. That is why the barrier condition is not needed for our uniqueness results. 
\end{remark}
{\bf Proof of Theorem \ref{MaxPrinc}:} 
We may assume that $E_i$ is open, $i=1,2$, since
otherwise we may replace 
$E_i$ by  $\mbox{int} \,E_i^{(1)}$, which
in view of Theorem \ref{singularity-estim} 
differs from $E_i^{(1)}$, and hence $E_i$, on a set of measure zero.
Also, if $E_i = \mbox{int} E_i^{(1)}$, then
clearly $\partial E_i = \partial E_i^{(1)}$,
and $E^{(1)}\subset E^{(2)}$ if  $E_1\subset E_2$.
So in the sequel we may  drop all superscripts on $E_i^{(1)}$, $i=1,2$ 
(but not for example on
$(E_1\cup E_2)^{(1)}$, since in  general $\partial (E_1^{(1)} \cup E_2^{(1)}) \ne \partial (E_1\cup E_2)^{(1)}$.)
 
It also suffices to prove the theorem under hypothesis \eqref{EmptyIntersec}, %eqref{EmptyIntersec},
since this follows from the barrier condition if $E_1 \setminus \Omega \subset \subset E_2 \setminus \Omega$.
Indeed, when (\ref{CC}) holds it is clear that $\partial E_1\setminus E_2 \subset \overline\Omega$. And if 
$x_0\in (\partial E_1 \setminus E_2) \cap \partial \Omega$,
then there exists $\epsilon_0>0$ such that $\partial E_1 \cap B(\epsilon, x_0) \subset \bar{\Omega}$ for all $\epsilon<\epsilon_0$, because otherwise $x_0 \in (\overline{E_1 \setminus \Omega}) \setminus E_2$ which violates the assumption $E_1 \setminus \Omega \subset \subset E_2 \setminus \Omega$.  On the other hand, if the barrier condition holds, then  according to Lemma \ref{SBC-Lemma},
it cannot be the case that
$\partial E_1 \cap B(\epsilon, x_0) \subset \bar{\Omega}$.
The proof of the other inclusion in \eqref{EmptyIntersec} is
essentially identical.

\medskip

\begin{figure}[ht]
   \begin{center}
	 \def \svgwidth{4cm}
	       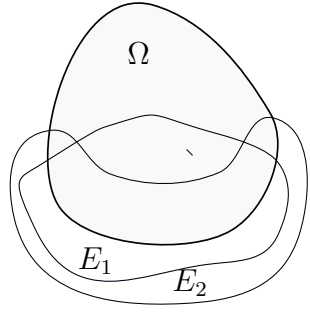
   \end{center}
   \caption{Failure of the comparison principle}
   \label{FWMaxPrinc}
\end{figure}

We now  assume toward a contradiction that   $E_1 \not \subset E_2$ (Figure \ref{FWMaxPrinc}). 
Note that Theorem \ref{singularity-estim} implies that $E_i \cap \Omega$
differs from its interior by a set of measure zero, so it follows that $E_1\setminus E_2$
has nonempty interior.
We claim that then 
\begin{equation}\label{Hdim}
dim_{Haus}(\partial E_1 \cap \partial E_2) \geq n-2.
%\calH^{n-2}(\partial E_1 \cap \partial E_2) >0.
\end{equation}
We will assume that $\calH^{n-1}(\partial E_1 \cap \partial E_2)=0$, since
otherwise \eqref{Hdim} is immediate.
Then either $\partial E_1 \setminus \overline E_2 $ or $\partial E_2\cap E_1$
must have positive $\calH^{n-1}$ measure,
since
$\partial(E_1\setminus E_2)$ is an $(n-1)$-dimensional set (by Proposition \ref{openProp}), and
\[
\partial (E_1\setminus  E_2) \subset (\partial E_1\setminus  \overline E_2) \cup (E_1\cap \partial E_2)
\cup (\partial E_1 \cap \partial E_2).
\]
For concreteness we assume that $\calH^{n-1}(\partial E_1 \setminus \overline E_2 ) >0$; the
other case is essentially the same.

Since $\partial E_1\setminus \overline E_2\subset \Omega$ by \eqref{EmptyIntersec}, 
and since Theorem \ref{singularity-estim} 
guarantees that  $\calH^{n-3}$ a.e. point of $\partial E_1\cap \Omega$ is regular,
%we find that $(reg(\partial E_1) \cap \Omega) \setminus \overline E_2$ is nonempty. In particular,
we may fix a connected component $R$ of $reg(\partial E_1)\cap \Omega$ such
that $R\setminus \overline E_2$ is nonempty.

It follows from Lemma \ref{intersectLem} that  $\bar R \cap \partial \Omega \neq \emptyset$, and
then  
$\bar R \cap \partial \Omega \subset \partial E_1\cap \partial \Omega \subset E_2$ by
\eqref{EmptyIntersec}.
Thus $\bar R \cap E_2 \ne \emptyset$. Since $E_2$ is open, it follows that $R\cap E_2\ne \emptyset$.

Then $R\cap \partial E_2$ separates $R$ into two nonempty components, $R\setminus \overline E_2$
and $R\cap E_2$.
Since the definition of a regular point implies that $R$ is an $(n-1)$-dimensional submanifold of $\R^{n}$,
and $R$ is connected by definition,
it follows from Proposition \ref{disconnect-prop} that the topological dimension of $R \cap \partial E_2$ is at least $n-2$. Since $R\cap \partial E_2\subset \partial E_1\cap \partial E_2$,
claim \eqref{Hdim} now follows from Proposition \ref{HTdim}.

It follows from  \eqref{Hdim} and Theorem \ref{singularity-estim}  that
\begin{equation}
\calH^{n-2} \big( reg(\partial E_1) \cap reg(\partial E_2)\big) >0.
\label{regreg}\end{equation}
Note also that
\begin{multline}
reg(\partial E_1) \cap reg(\partial E_2) \ \subset \ 
\left[ \partial(E_1\cup E_2)^{(1)} \cap \partial (E_1\cap E_2)^{(1)}
\right]\\
\cup
\left[ \partial(E_1\cup (\overline E_2)^c)^{(1)} \cap \partial (E_1\cap (\overline E_2)^c)^{(1)}
\right].
\label{cupcap}\end{multline}
Indeed, if $x_0 \in  reg(\partial E_1) \cap reg(\partial E_2)$,
then by the definition of a regular point,
we may write both boundaries $\partial E_1$ and
$\partial E_2$ 
locally as graphs
over the same domain in the same hyperplane, say of $C^2$ functions
$w_1$ and $w_2$. We may also assume that $E_1$ lies below the graph
of $w_1$ near $x_0$. It is straightforward to verify that if $E_2$ also lies below the graph
of $w_2$ near $x_0$, then
\[
x_0\in  \partial(E_1\cup E_2)^{(1)} \cap \partial (E_1\cap E_2)^{(1)}
\]
whereas if $E_2$ lies above the graph of $w_2$ near $x_0$, then
\[
x_0 \in \partial(E_1\cup (\overline E_2)^c)^{(1)} \cap \partial (E_1\cap (\overline E_2)^c)^{(1)}.
\]

By Corollary \ref{areaMinCor}, all of the sets $E_1\cup E_2$, $E_1\cap E_2$, \ldots on the right-hand side of \eqref{cupcap}
are $\varphi$-area minimizing on $\Omega$, and hence their boundaries are
all regular in $\Omega$ away from a set of dimension at most $n-3$.
It follows that for a suitable choice $F_1 = E_2$ or $\bar E_2^c$, 
\eqref{Hdim} and \eqref{DimSingSet} that
\begin{equation}
U:=\{x\in \partial E_1 \cap \partial E_2: \ \ \partial (E_1 \cup F_1)^{(1)}  \ \ \hbox{and} \ \ \partial (E_1 \cap F_1)^{(1)} \ \ \hbox{are regular at} \ \ x\}
\label{U.def}\end{equation}
satisfies $\calH^{n-2}(U)>0$.
Let $R$ be a connected component of $reg(\partial E_1)$ that intersects $U$, and let 
\[
R_0 := R \setminus S, \quad \qquad
S := 
 sing(  \partial (E_1 \cup F_1)^{(1)}) 
\cup 
 sing(  \partial (E_1 \cap F_1)^{(1)}) .
\]
We claim that $R_0\subset U$.
To prove this, first note that 
the topological dimension of $S$ is bounded by $dim_{Haus}(S) \le n-3$, 
so Proposition \ref{disconnect-prop} implies that $R_0$ is connected.
Thus it suffices to show that $U\cap R_0$ is nonempty,  open and closed in $R_0$.
It follows from the definitions $\emptyset \ne U\cap R \subset U\cap R_0$,
and also that $U\cap R_0 = R_0 \cap \partial E_2$, which is
clearly closed.
So we only need to check openness. For this, fix $x\in R_0\cap U$,
and note that
since $(E_1 \cap F_1) \subset  (E_1 \cup F_1)$
and both boundaries are
regular at $x$, we may write both boundaries locally as $C^2$ graphs
over the same domain in the same hyperplane, say of functions
$w_1$ and $w_2$
such that $w_1\le w_2$ in their domain. Then 
Lemma \ref{HopfMax} implies that  $\partial (E_1\cap F_1)$ coincides with $\partial (E_1 \cup F_1)$ in a neighbourhood of $x$, which implies that $E_1$ coincides with $F_1$
in a neighborhood of $x$ and hence that $\partial E_1$ and $\partial E_2 = \partial F_1$
coincide in a neighborhood of $x$. Thus $U\cap R_0$ is open.
 
Now the dimension estimate of $S$ implies that $R_0$ is dense in $R$,
and thus $\overline R  = \overline R_0 \subset \overline U$.
We deduce using Lemma \ref{intersectLem},
that $ \partial E_1\cap \partial E_2 \cap \partial \Omega$ is nonempty,
which is impossible due to \eqref{EmptyIntersec}. Thus we have arrived at a conradiction.

%\begin{eqnarray*} %\label{DimSingSet}
%\left\{ \begin{array}{ll}
%\mathcal{H}^{n-3}(S)< \infty \ \ &\hbox{if} \ \ n\geq 4\\
%S =\emptyset, \ \ &\hbox{if} \ \ n\le 3.
%\end{array} \right.
%\end{eqnarray*}

Finally, suppose that $n\le 3$. We already know that $E_1\subset E_2$, and
if $E_1$ is not compactly contained
in $E_2$, then $\partial E_1 \cap \partial E_2\cap \Omega$ is nonempty. 
Since $n\le 3$, $\partial E_1$ and $\partial E_2$ are both regular everywhere in
$\Omega$, so we can  invoke Lemma \ref{HopfMax} to find that
$U := \partial E_1\cap \partial E_2$ is open in both  $\partial E_1$ and $\partial E_2$.
Then arguing as above, we find that $\partial U \cap \Omega$ is nonempty, and this again is a contradiction.
\hfill $\Box$

\ \ 

Now we use Theorem \ref{MaxPrinc} to establish our main
uniqueness and continuity results.\\ \\
{\bf Proof of Theorem \ref{Comparison}:}  We first prove (\ref{comparison}). As before we extend $f_i$ for $i=1,2$ to continuous functions on $\overline{\Omega^c}$ (still denoted by $f_i$). 
We may assume that $f_1 \le f_2$ in $\Omega^c$, since otherwise we may
replace $f_1, f_2$ by $\min(f_1, f_2)$ and $\max(f_1,f_2)$ respectively. 
For $i=1,2$, we extend $u_i$ to a function  (still denoted $u_i$) 
on $\R^n$ by setting it
equal to $f_i$ on $\Omega^c$.

Suppose toward a contradiction that
\eqref{comparison} is not true.
Then since
\[
\{x\in \Omega : u_1(x) > u_2(x) \} = \bigcup_{(\lambda_1, \lambda_2)\in \Q\times \Q   }
\left\{ x\in \Omega : u_1(x) \ge \lambda_1 > \lambda_2 > u_2(x)\right\}
\]
there must be some rational numbers $\lambda_1 > \lambda_2$ such that
\[
\calH^n\left( \left\{ x\in \Omega : u_1(x) \ge \lambda_1 > \lambda_2 > u_2(x)\right\} \right) >0.
\]
If we define
\[
E_i := \{ x\in \R^n : u_i(x) \ge \lambda_i\}.
\]
then this says exactly that $\calH^n(E_1 \setminus E_2) >0$. 

On the other hand, it follows from the definition of $BV_{f_i}(\Omega)$
and the continuity of $f$
that if  $x\in \partial E_i^{(1)}\cap \partial \Omega$,
then $f_i(x) = \lambda_i$. Since $f_1 \le f_2$ on $\partial \Omega$ and $\lambda_1>\lambda_2$,
the definitions imply that \eqref{EmptyIntersec} holds.
Hence it follows from Theorem \ref{MaxPrinc}  that $E_1^{(1)}\subset E_2^{(1)}$,
and therefore that 
$\calH^n (E_1 \setminus E_2)=0$. Since this contradicts the above,
we conclude that \eqref{comparison} holds.

Finally, as is well known, it is easy to deduce 
\eqref{uniqueness} from \eqref{comparison}.
For example, to prove that
$u_1 - u_2 \le  \sup_{\partial \Omega} |f_2-f_1|$ {\em a.e. } in $\Omega$, we apply
\eqref{comparison} to $u_1$ and $\tilde u_2 := u_2 +\sup_{\partial \Omega} |f_2-f_1|$.
The opposite inequality is proved by the same argument.
\hfill $\Box$

\begin{remark}\label{X}
For $\varphi(x,\xi) = a(x)|\xi|$, or more generally $\varphi$ of the
form (\ref{WLG}), minimizers are continuous in $\Omega \subset R^n$ for $n \le 7$.
The point is that in this case the boundary of a  $\varphi$-area
minimizing set is actually a minimal hypersurfaces with  respect to some Riemannian metric, and
as such has better regularity properties than in the case of the more general class of integrands we consider in this paper. In
particular, if $E$ is $\varphi-$area minimizing in $\Omega\subset R^n$ and $n \le 7$, then
$sing(\partial E^{(1)}) \cap \Omega = \emptyset$. (This is documented for example in \cite{FMorgan}.) Continuity for $n \le 7$ is established by using this fact in place of
Theorem \ref{singularity-estim} and repeating the above arguments.
\end{remark}

\vspace{0.5cm}
{\bf Proof of Theorem \ref{thm.cont}:}
As in Theorem \ref{Comparison}, we set $u$ equal to $f$ on $\Omega^c$.
For $x\in \Omega$ we define
\[
u^* (x):= \lim_{r\to 0} \esssup_{B(r,x)}u,
\qquad
u_* (x):= \lim_{r\to 0} \essinf_{B(r,x)}u,
\]
We must show that $u^*= u_*$ everywhere in $\Omega$.
Assume toward a contradiction that this fails, so that
\[
u_*(x_0) < \lambda_2  < \lambda_1 < u^*(x_0)
\]
for some $x_0\in \Omega$ and $\lambda_1,\lambda_2\in\R$.
Define 
\[
E_i := \{ x\in \R^n : u(x) \ge \lambda_i \}
\]
As in the proof of Theorem \ref{Comparison},
the hypotheses of Theorem \ref{MaxPrinc}
are satisfied, and (since now $n\le 3$) 
it follows that $E_1^{(1)} \subset\subset E_2^{(1)}$,
and hence that $\partial E_1^{(1)} \cap \partial E_2^{(1)} = \emptyset$.

So to arrive at a contradiction, it suffices to check that $x_0\in \partial E_i^{(1)}$ for $i=1,2$.
This is straightforward. In fact, since $u^*(x_0) > \lambda_i$, every ball around $x_0$ contains
a subset of $E_i$ of positive measure, and hence a subset of $E_i^{(1)}$ of (the same) positive
measure.
Similarly, since $u_*(x_0) < \lambda_i$, every ball around $x_0$ contains a
subset of $(E_i^{1})^c$ of positive measure, and it follows that $x_0\in \partial E_i^{(1)}$ for $i=1,2$.
\hfill $\Box$

\medskip

Finally, we conclude this section with the proof of Lemma \ref{intersectLem},
which played an important role in the above arguments.\\ \\ 
{\bf Proof of Lemma \ref{intersectLem}}:
Let $R$ be a nonempty connected component of $reg (\partial E^{(1)}) \cap \Omega$,
for a set $E$ that is $\varphi$ area-minimizing in $\Omega$.
Also, assume toward a contradiction that $\bar R \cap \partial \Omega = \emptyset$.

The definition of a regular point implies that $R$ is a $C^2$ submanifold
of $\R^n$.
Given a smooth $(n-1)$-form $\psi$ in $\R^n$, we will write 
\[
\bbl R \bbr (\psi) := \int_R \psi
\]
where $R$ has the same orientation as $\partial E$, which in turn inherits
its orientation from $E$ in the standard way.
And given an $(n-2)$-form $\eta$, we will write
\begin{equation}
\partial \bbl R\bbr(\eta) := \bbl R \bbr(d\eta) .
\label{dR}\end{equation}
(In particular, here ``$\partial$" denotes the boundary in a distributional sense, rather than
the topological boundary.) 
%In the language of geometric measure theory, $\bbl R \bbr$ and $\partial \bbl R \bbr$ arean $(n-1)$-current and an $(n-2)$-current respectively. 
We first claim that 
\begin{equation}
\partial \bbl R \bbr  = 0.
\label{dRzero}\end{equation}
Toward this end,
note that Stokes' Theorem implies that
$\partial \bbl R \bbr (\eta ) = 0$ if the support of $\eta$ does not intersect $\overline R \setminus R$.
Since $\overline R  \subset \Omega$, the definitions imply that
$\overline R\setminus R \subset \partial E \setminus reg(\partial E) = sing(\partial E)  =: S$.
Thus $\partial \bbl R \bbr$ is supported in $S$. Since $\calH^{n-2}(S) =0$,
the claim \eqref{dRzero} follows from standard geometric measure theory considerations,
which we summarize as follows:
\begin{itemize}
\item It is clear that $\bbl R \bbr$ is a $(n-1)$-dimensional rectifiable current and
hence a $(n-1)$-dimensional flat chain; see \cite[4.1.24]{federer}.
\item Thus $\partial \bbl R \bbr$ is a $(n-2)$-dimensional flat chain, see \cite[4.1.12]{federer}.
\item As an  $(n-2)$-dimensional flat chain whose support has $\calH^{n-2}$ measure zero,
$\partial \bbl R \bbr$
must be trivial, see \cite[4.1.20 and 2.10.6]{federer}. This is \eqref{dRzero}.
\end{itemize}

Next we claim that there exists a set $F\subset\subset\Omega$ of finite perimeter
such that 
\begin{equation}
\bbl R \bbr \ = \ \pm \partial \bbl F \bbr
\label{R.bdy2}\end{equation}
where, on the right-hand side, $\bbl F \bbr$ denotes the $n$-current corresponding to
integration over $F$, and $\partial$ is defined as in \eqref{dR}. (The meaning of $\bbl \,\cdot\, \bbr$ should always be clear from the context.)
Indeed, the fact that $R$ is bounded and $\partial \bbl R \bbr = 0$ implies
that there exists some compactly supported integer-valued BV function $u_R$
such that 
\begin{equation}
\bbl R \bbr( \psi)  \  =  \ \int_{\R^n} u_R \  d\psi
\qquad
\mbox{ for every smooth compactly supported $(n-1)$-form $\psi$.}
\label{uF.def}\end{equation}
This follows for example as an easy special case of the isoperimetric theorem, see \cite[4.2.10]{federer}.
One can deduce from
\eqref{uF.def} and the coarea formula (or for a detailed proof see see \cite[Theorem 27.6]{simon}) that
\begin{equation}
\bbl R \bbr = \sum_{k=1}^\infty \partial\bbl  \{ x : u_{R}(x) \ge k \} \bbr - \sum_{k=1}^\infty \partial
\bbl  \{ x : u_{R}(x) \ge -k\} \bbr
\label{decomp1}\end{equation}
and 
\begin{equation}
\calH^{n-1}(R) = \sum_{k=1}^\infty \calH^{n-1}( \partial \{ x : u_{R}(x) \ge k \} ) +
 \sum_{k=1}^\infty  \calH^{n-1}( \partial \{ x : u_R{}(x) \le -k\}).
\label{decomp2}\end{equation}
However, it is shown in \cite[4.2.25]{federer}
that the fact that $R$ is a $C^2$, connected submanifold
of $\R^n$ 
implies that 
there can only be one nontrivial term on the right-hand side of 
\eqref{decomp1}. %(The claim there is if $R$ is a connected $C^1$ submanifold  of $\R^n$, then $\bbl R \bbr$ is indecomposable, and by definition, this means that it admits no nontrivial decomposition \eqref{decomp1} satisfying \eqref{decomp2}.) 
(The point is that a smooth connected  submanifold of $\R^n$ is {\em indecomposable},
which means exactly that it admits no nontrivial decomposition \eqref{decomp1} satisfying \eqref{decomp2}.) 
So $u_R$ is the characteristic function of a set $F$,
up to a sign, and \eqref{uF.def} thus reduces to the claim  \eqref{R.bdy2}. We must also show that $F\subset\subset\Omega$.
To see this, we observe from
\eqref{uF.def} that
$Du_R=0$  in the sense 
of distributions
away from $R$, and in particular in $\Omega^c$. The fact that $\partial \Omega$ is connected implies 
that  $\Omega^c$ is connected.
Hence, since it has compact support, $u_R = 0$ in $\Omega^c$, so
that $F \subset \bar \Omega$.
Moreover, since $\overline R\cap \partial \Omega = \emptyset$
and $\overline R, \partial \Omega$ are both compact, these sets are separated
by a positive distance. Hence, reasoning as above we see that $u_R=0$
in a neighborhood of $\partial \Omega$, which means that $F\subset\subset \Omega$.

By duality between $(n-1)$-forms and vector fields, \eqref{R.bdy2} is equivalent to
\[
\int_R \eta\cdot \nu_{E} \  d\calH^{n-1} \ = \ 
\pm \int_F \nabla \cdot \eta \  d\calL^{n} 
\qquad
\mbox{ for all } \eta\in C^\infty_c(\R^n;\R^n),
\]
and in view of the smoothness of $R$, this implies that 
\begin{equation}
\overline R  = \partial F^{(1)}, \qquad \nu_{E} = \pm \nu_F \mbox{ \ on \ }R.
\label{RfD}\end{equation}
First assume for concreteness that $ \nu_{E} = \nu_F$ on $R$, and define
$\wt E := E\setminus F$.
Then $E \cap \Omega^c = \wt E\cap \Omega^c$. We will show that
$P_\varphi(\wt E)< P_\varphi(E)$, contradicting the minimality of $E$.
We will use the notation 
\[
E^{(0)}:= \{x\in \R^n :
\lim_{r\to 0} \frac {\calH^n(B(r,x)\cap E)}{\calH^n(B(r))} = 0 \},
%\qquad \qquad E^{(0)} := (E^c)^{(1)}.
\qquad
\partial_* E := \R^n \setminus (E^{(0)}  \cup E^{(1)})
\]
and the fact that
for any set $E$  of finite perimeter,
$\partial_* E$ is $\calH^{n-1}$ measurable, an approximate unit
normal $\nu_E$ exists $\calH^{n-1}$ a.e. in $\partial_*E$,  and
\[
P_\varphi(E) := \int_{\partial_*E} \varphi(x,\nu_E) d\calH^{n-1}.
\]
For the sets $E,F$ above, the regularity of $R$ implies that $\partial_* E =
\partial E^{(1)}$ in $\Omega$, up to sets of dimension $n-3$, and similarly for $F$.
In addition,
\[
\partial_*\wt E \ = \ 
\partial_*(E\setminus F) \ = \ 
(\partial_*E \cap F^{(0)}) \  \cup \  (
E^{(1)}\cap \partial_* F) \ \cup \ \{ x\in  \partial_*E\cap \partial_* F : \nu_{E} = - \nu_{F}\}
\]
up to sets of $\calH^{n-1}$ measure zero. 
The content of this statement can
be understood by drawing a picture, and a proof can be found for example in
\cite{Maggi}, Theorem 16.3. Then since $\nu_E = \nu_F$ on $R$,
it follows that 
\begin{equation}
\partial_* \wt E  \ = \  \partial_* E\cap F^{(0)} \  \subset \ 
\partial_*  E \setminus R \  \subset \ 
 \partial_* E,
% \subset \ \calH^{n-1}(\partial_*  E \setminus \partial_*  \wt E)  >0
\label{EminusE}\end{equation}
up to $\calH^{n-1}$ null sets.
Then we complete the proof by calculating
\[
P_\varphi(E) 
= 
\int_{\partial_*E} \varphi(x,\nu_{E}) d\calH^{n-1}
>
\int_{\partial_*E\setminus R} \varphi(x,\nu_{E}) d\calH^{n-1}
\ge
\int_{\partial_*\wt E} \varphi(x,\nu_{\wt E}) d\calH^{n-1}
=
P_\varphi(\wt E),
\]
If $\nu_E = -\nu_F$ on $R$ then we define $\wt E := E\cup F$
and argue in essentially the same way.

%%%%%%%%%
%%%%%%%%%
\hfill $\Box$

\section{Non-uniqueness}\label{S:counter}

In this section we show that the regularity  assumptions in our uniqueness theorems are sharp.
Indeed for any $\alpha<1$ and $n\geq 2$, we will construct a $ C^{1,\alpha}$ function $a$ on a bounded  region $\Omega \subset \R^n$ such that  the weighted least gradient problem (\ref{LTVProb}) has infinitely many minimizers in $BV(\Omega)$ (Proposition \ref{non-unique-prop}).

The main ingredient in the proof is the following Lemma,
in which we simultaneously construct both a family of functions $u_\sigma$
and a vector field that calibrates them.
In the lemma, we will write points in $ \R^n$ in the form $(x,z)\in \R^{n-1}\times \R$.

\begin{lemma}
Let $D := \{(x,z)\in \R^{n-1}\times \R : |x| < 3 \}$. For $0\le \sigma \le 1$,
there exist a function $u_\sigma\in BV(D)$
and a vector field $J\in C^{\frac {1+\alpha}2}(D;\R^n)$
such that
\begin{equation}\label{J1}
\nabla \cdot J = 0 \quad\mbox{ in }\calD'(D),
\qquad 
\qquad |J| \in C^{1,\alpha}(D),
\end{equation}
\begin{equation}\label{calib-ident}
\int_{D}J\cdot Du_{\sigma}= \int_{D}|J| \ |Du_{\sigma}|,
\end{equation}
and $|J|>0$ in $D$. %$ is bounded away from zero on a neighbourhood of $supp(Du_{\sigma})$, for all $\sigma\in [0,1]$. 
\label{calib-lemma}
\end{lemma}

We will need not just that $\nabla\cdot J$ vanishes in the sense of distributions, but
actually
a slightly stronger property. This will be established later, in the proof of the Proposition,
see \eqref{IBP}, so we do not carefully check this condition here.

\medskip

%\begin{lemma}\label{calib-lemma}
%For any $\alpha<1$ there exists a vector field 
%$J\in C^{\frac{1+\alpha}{2}}(D)$ with $|J| \in C^{1,\alpha}(D)$ 
%such that $\nabla \cdot J\equiv0$,
%\begin{equation}\label{calib-ident}
%\int_{D}J\cdot Du_{\sigma}= \int_{D}|J||Du_{\sigma}|,
%\end{equation}
%and $|J|$ is bounded away from zero on a neighbourhood of $supp(Du_{\sigma})$ for all $\sigma \in \R$. 
%\end{lemma}

{\bf Proof: }  
For $(x,z)\in D$, we will write $r := |x|$. We will take
$J$ to have the form
\begin{equation}\label{J}
J(x,z)=(\frac{-x}{|x|^{n-1}} \psi_z(x,z ), \frac{1}{|x|^{n-2} }\psi_r (|x|,z)),
\end{equation}
for a function $\psi(r,z)$ to be chosen,
satisfying
\begin{equation}
\psi(r,z)= r^{n-1} \quad\mbox{ if }r \le \frac 18
\label{psi1}\end{equation}
as well as other conditions to be stated later.
Such vector fields always satisfy
$\nabla \cdot J = 0$ in the sense of distributions,
provided $\psi$ is regular enough, which will be the
case here. (See \eqref{IBP} below.)
So the main point is now simply 
\eqref{calib-ident}.  To arrange that this holds,
we will also take $u$ to depend only on $r$ and $z$,
in which case
\eqref{calib-ident} reduces to the condition 
that level curves of $u_\sigma$ in the $r\!\!-\!\!z$ plane are 
orthogonal to $(-\psi_z, \psi_r)$, or equivalently, are parallel to $\nabla \psi$. 
This is an
ordinary differential equation for level curves of $u_\sigma$,
and if $\nabla\psi$ is not Lipschitz continuous, solutions
need not be unique. Thus failure of uniqueness on the level
of ODEs will be the basis for failure of uniqueness 
for the variational problem.

Our first task is thus to choose $\psi$ so that uniqueness
fails for certain flows along the vector field $\nabla\psi$.
We find it easier to write down a
(mostly) explicit example, chosen
to facilitate our later construction of $u_\sigma$,
than to proceed by abstract arguments.
So in addition to \eqref{psi1},
we require that
\begin{equation}
\psi(r,z) =
r - g(r)(\frac 1{1+\theta})|z|^{1+\theta} \quad \mbox{ if }r \ge \frac 14
\label{psi2}\end{equation}
for $\theta<1$ to be fixed below, where $g$ is a smooth function such that 
\[ %\begin{equation}\label{u}
g(r)=
\begin{cases}
0 \ \ &\hbox{if} \ \ r \le \frac 13, \\
1\ \ &\hbox{if} \ \ r\in [\frac{1}{2},1],\\
-1\ \ &\hbox{if} \ \ r\in [2,3],
\end{cases}
\qquad\qquad  \mbox{ and } \ \ 
g(r) \ge 0 \ \ \hbox{ if} \ \ r\in (\frac 13, \frac 12).
\] %\end{equation} 
We finally require that $\psi(r,z) $ is a smooth function of $r$ alone when $r\le \frac 14$,
such that \eqref{psi1} holds and $\psi_r>0$ for $0<r\le \frac 14$.

It is clear that $J$ is smooth away from $\{ z=0, r\ge \frac 14\}$,
and also that $J$ is $C^{0,\theta}$.
Also, for $r \ge \frac 13$,
\[
|J| =  r^{2-2n}|\nabla\psi| =  r^{2-2n}( \psi_r^2 + g^2(r) z^{2\theta})^{1/2}
\]
and since $\psi_r^2$ is $C^2$ and positive, it is straightforward  to check that
$|J| \in C^{1, 2\theta -1}(D)$. To satisfy the conclusions of the lemma we choose $\theta=\frac{1+\alpha}{2}$.

Now we define a family of functions
$u_\sigma(r,z)$ such that every level curve
of every $u_\sigma$ 
is an integral curve of $\nabla\psi$.

\begin{figure}[ht]
   \begin{center}
	 \def \svgwidth{10cm}
    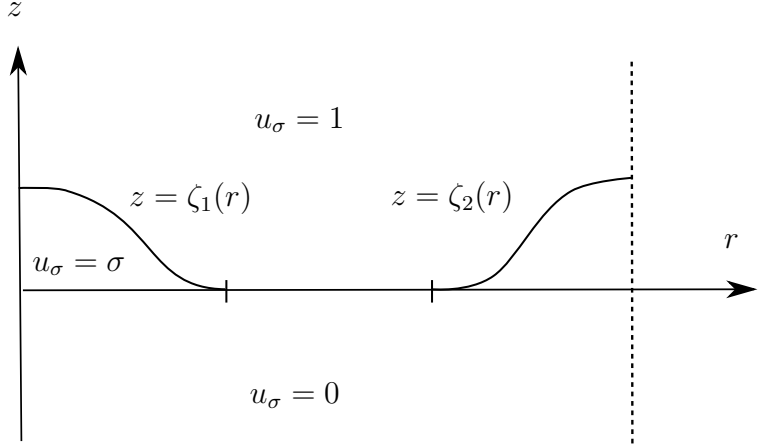
   \end{center}
   \caption{$u_\sigma(r,z)$}
   \label{CE}
\end{figure}

Toward this end, note that $\nabla\psi = (\psi_r, 0)$ when $z=0$, with $\psi_r >0$,
so the $r$-axis is an integral curve of $\nabla \psi$.
Notice also that $\nabla\psi = (1,  \pm z^\theta)$ if $r\in [\frac 12,1]\cup  [2,3]$,
and so in these regions one can explicitly integrate to find
integral curves of $\nabla \psi$.
The following definition (in which 
$\zeta_1,\zeta_2$ will be defined below)
is thus quite natural:
\begin{align*}
u_\sigma(r,z) &:= 0  \hspace{10em}\mbox{ if }z<0 , \\
u_\sigma(r,z) &= 1 \hspace{10.5em} \mbox{ if }
\begin{cases}
r \le 1 \mbox{ and }z> \zeta_1(r) \quad &\\
1 \le r \le 2 \mbox{ and }z>0, \quad &\\
r \ge 2 \mbox{ and }z> \zeta_2(r) \quad &
\end{cases} \\
u_\sigma(r,z) &:=\sigma
\hspace{10em}\mbox{ if }
r \le 1 \mbox{ and }0 < z < \zeta_1(r), \\
u_\sigma(r,z)&:=
\frac 1{1-\theta}z^{1-\theta}  - r + 3 
\hspace{3em}\mbox{ if }
r \ge 2 \mbox{ and }0 < z < \zeta_2(r).
 \end{align*}
Figure \ref{CE} is a sketch of these regions in the (right half of the) $r-z$ plane.

To complete the definition, we need to specify $\zeta_1(r), \zeta_2(r)$.
First, we  choose $\zeta_2$ such that $u_\sigma$
is continuous on the set $\{ (r, \zeta_2(r) : 2\le r \le 3\}$. 
Thus $\zeta_2(r)$ is characterized by the identity
\begin{equation}
\{ (r, \zeta_2(r)): 2\le r \le 3 \} := \{ (r,z) : 2\le r \le 3, 
\frac 1{1-\theta}z^{1-\theta}  - r + 3 
=1 \}.
\label{c2.def}\end{equation}
Second, we choose $\zeta_1(r)$ to be any function such that
$r\in [0,1] \mapsto (r, \zeta_1(r))$ is an integral curve of $\nabla\psi$
such that $\zeta_1(1) = 0$ and $\zeta_1>0$ somewhere in $(0,1)$.
For example, we choose $\zeta_1$
to equal $\left[ (1-\theta)(r-1)\right]^{1/(1-\theta)}$
for $\frac 12 \le r \le 1$. Then $\zeta_2(r)$ for 
$0\le r \le \frac 12$ may be found by solving the 
appropriate ODE, which is just $(1, \zeta_1'(r))\cdot \nabla^\perp\psi(r, \zeta_1(r))=0$.

It remains to verify that $u_\sigma\in BV(D)$ 
and that \eqref{calib-ident} holds for every $\sigma\in [0,1]$.
Let us write $\nabla u_{\sigma}$ to denote the absolutely continuous part of $Du_{\sigma}$ and $\Sigma$ to denote the jump set of $u_\sigma$, which is independent of $\sigma$ for 
$\sigma\in (0,1)$. Recalling that $r = |x|$ for $x\in\R^{n-1}$, we have
\[ %\begin{equation}\label{u}
\nabla u_{\sigma}=
\begin{cases}
( - \frac x r, z^{-\theta}) \  \ \ &\mbox{ if }
r \ge 2 \mbox{ and }0 < z < \zeta_2(r)\\
0&\mbox{ otherwise}
\end{cases}
\] %\end{equation}
and $\Sigma \subset \Sigma_1 \cup \Sigma_2$, where
\begin{equation}
\Sigma_1= \{(x,z)\in D : z = 0 \}, \qquad
\quad
\Sigma_2 := \{(x,z)\in D : r\le 1, z = \zeta_1(r)\}.
\label{Sigma.def}\end{equation}
(Indeed, $\Sigma = \Sigma_1$ if $\sigma=1$ and $\Sigma_2$ if $\sigma=0$, and
otherwise $\Sigma = \Sigma_1\cup \Sigma_2$.)
It is easy to check that $\nabla u_\sigma\in L^1$ and since $u_\sigma \in L^\infty$
and the jump set has finite $\calH^{n-1}$ measure, we infer that
$u_\sigma \in BV(D)$ for every $\sigma$, and
\begin{equation}\label{splitgrad}
\int_{D}J\cdot Du_{\sigma}=\int_{D} J\cdot \nabla u_{\sigma}+\int_{\Sigma} J\cdot \nu[u_{\sigma}^+ - u_{\sigma}^-]d \calH^{n-1},
\end{equation}
where $\nu$ is the upward unit normal to $\Sigma$ and $u_{\sigma}^\pm(x,y)=\lim_{\epsilon \rightarrow 0^\pm}u_{\sigma}(x,y+\epsilon)$. 
We claim that
%It is easy to check  from the definitions that 
\begin{equation}\label{1}
J  \cdot \nabla u_{\sigma}=|J| |\nabla u_{\sigma}|\ \ \mbox{ in }D,
\qquad\qquad
J  \cdot \nu =|J  | \ \ \hbox{on} \ \ \Sigma.
\end{equation}
Indeed, the first claim is a routine verification, as is the second
for $\Sigma_1$. For $\Sigma_2$ it follows in a straightforward way 
from the ODE solved by $\zeta_1(r)$.

Finally, since $u_\sigma^+- u_\sigma^- \ge 0$ on $\Sigma$ for every $\sigma\in [0,1]$,
we deduce from \eqref{splitgrad} and \eqref{1} that \eqref{calib-ident} holds.
\hfill $\Box$

\medskip

We now show that the functions $(u_\sigma)_{\sigma\in [0,1]}$ constructed above
yield a counterexample to uniqueness for a suitable weight $a$.

\medskip 
 {\bf Proof of Proposition \ref{non-unique-prop}:} 
Let $\Omega$ be an open subset of $D$ that contains the jump set $\Sigma = \Sigma_1 \cup \Sigma_2$,
see \eqref{Sigma.def}. Then $u_\sigma$ is continuous at $\partial \Omega$ and $u_{\sigma}|_{\partial \Omega}=u_{0}|_{\partial \Omega}$ for all $\sigma \in [0,1]$. In particular, if we write $f$
to denote the boundary value of $u_0$, then $u_\sigma\in BV_f(\Omega)$ for
every $\sigma\in [0,1]$.

Next, let $a=|J|$ where $J$ is the vector field constructed in Lemma \ref{calib-lemma}.  We now claim that
that for every $w \in BV(\Omega)$, we have 
 \begin{eqnarray}\label{IBP}
\int_{\Omega} J\cdot Dw &=& \int_{\partial \Omega} w\,J\cdot \nu.
\end{eqnarray}
To see this, write $J=(J_1,J_2)$ with $J_1$ denoting the first $n-1$ components of $J$.
%; note that $J_1\in C^{1}(\{(x,z)\in \Omega  : z=0\})$, $\partial_x J_1\in C(\Omega)$, and $J_2\in C^{1}(\Omega)$.  
Let $\rho:\R\rightarrow \R$ be a smooth function such that $\rho\equiv 0$ in $[-1,1]$ and $\rho = 1$ for $|t|\geq 2$. Define $\rho_\epsilon(t)=\rho(\frac{t}{\epsilon})$ and 
\[J_\epsilon(x,z)=(J_1(x,z) \rho_\epsilon(z), J_2(x,z)).\]
Then from the  form \eqref{J} of $J$ one checks that $J_\epsilon \in C^{1}(\Omega)$ and $\nabla \cdot J_\epsilon =0$, so that
\[
\int_{\Omega} J_\epsilon\cdot Dw \ = \ \int_{\partial \Omega} wJ_\epsilon\cdot \nu.
\qquad\mbox{ for every }w\in BV(\Omega).
\]
Since $\rho_\epsilon$ is independent of $x$, letting $\epsilon \rightarrow 0$ we obtain (\ref{IBP}).
It follows that
\[
\int_{\Omega} a|Du_\sigma| \ \overset{\eqref{calib-ident}} = \  \int_{\Omega} J\cdot Du_\sigma
\ = \  \int_{\partial \Omega}J\cdot \nu u_0.
\]
Now consider any $w\in BV_f(\Omega)$. Then $w|_{\partial \Omega} = f = u_0|_{\partial \Omega}$,
so by (\ref{IBP}) we obtain 
\[
 \int_{\Omega} a|Dw| 
\ \ge \  \int_{\Omega} J\cdot Dw
\ = \  \int_{\partial \Omega}J\cdot \nu\,  u_0
\ = \  \int_{\Omega} a|Du_\sigma|
\]
for all $\sigma \in [0,1]$. Therefore for every $\sigma \in [0,1] $, the function $u_\sigma$ is a minimizer of weighted the least gradient problem \eqref{WLG}.
%
%\begin{equation*}\label{TVProblem2}
%\inf  \left\{\int_{\Omega} a|Dw|: \ \ w\in BV_f(\Omega) \right \}.
%\end{equation*}
\hfill $\Box$

\section{About the barrier condition}\label{S:barrier}

This section is devoted the proof of the lemma, stated earlier, that 
provides a reformulation of the barrier condition for smooth domains.
%We first recall the following standard consequence of strict convexity.

%\begin{lemma} Let $H^- := \{x\in \R^n : x_n<0\}$, and assume that $\varphi_0:\R^n\to\R^n$ is a norm on $\R^n$ that is strictly convex in the sense that if $\xi\in S^{n-1}$, then
%\begin{equation}
%\xi \cdot D\varphi_0(\nu)  \le \varphi_0(\xi), \quad\mbox{ with equality if and only if }\xi = \nu
%\label{sconvex}\end{equation}.
%Let $Z\subset \R^n$ satisfy $Z\setminus B(1) = H^-\setminus B(1)$. Then
%\[
%P_{\varphi_0}(Z; B(2)) \ge P_{\varphi_0}(H^-, B(2)),
%\]
%with equality iff $Z = H^-$ a.e. 
%\label{L.convex.min}\end{lemma}

%{\bf Proof}
%We write $\zeta := D\varphi(e_n)$. Then, writing $\nu_Z $ for the Radon-Nikodym derivative $\frac{d\, D\chi_Z}{d\, |D\chi_Z|}$,  we estimate
%\[ 
%P_{\varphi_0}(Z;B(2)) =  \int_{B(2)} \varphi_0(\nu_Z) | D \chi_Z|  \ge  \int_{B(2)} \zeta\cdot \nu_Z  | D \chi_Z| = \int_{B_2} \zeta \cdot D\chi_Z.
%\] 
%But since $Z=H^-$ a.e. outside $B(1)$ and $\zeta\cdot \nu_{H^-} = \zeta\cdot e_n = \varphi_0(e_n)$ everywhere  in $\partial H^-$,
%\[ 
%\int_{B_2} \zeta \cdot D\chi_Z = \int_{B_2} \zeta \cdot D\chi_{H^-} = \int_{B(2)} \varphi(\nu_{H^-}) \, |D\chi_{H^-}| = P_{\varphi_0}(H^-;B(2)) .
%\]
%By inspecting the above argument and using again \eqref{sconvex} and the fact that $Z$ and $H^-$ coincide outside $B(1)$, one can check that the inequality is strict if and only if $Z = H^-$ a.e..
%\hfill $\Box$

\medskip

{\bf Proof of Lemma \ref{L.barrier}:}  
Assume toward a contradiction that \eqref{barrier.suff} holds, but that
the barrier condition fails. Then there exists 
$x_0\in \partial \Omega$, a sequence $\e_k\to 0$, and  sets
$V_k \subset \Omega$
such that
\[
V_k \mbox{ minimizes  }P_\varphi(\cdot, \R^n) \ \mbox{ in } \ 
\left\{ W\subset\Omega : W \setminus B(\e_k,x_0)  = \Omega \setminus B(\e_k,x_0) \right\},
\]
but
\begin{equation}
\partial V_k^{(1)} \cap \partial \Omega \cap B(\epsilon_k, x_0) \ne \emptyset.
\label{barrier.fail}\end{equation}
We will replace $V_k$ by $V_k^{(1)}$ and drop the superscripts.
By a change of coordinates, we may assume that $x_0=0$, and that 
$\nu(x_0)= (0,\ldots, 0,1)$, where $\nu$ denotes the outer unit normal
to $\partial \Omega$.
For each $k$ we define
\[
Z_k := \frac 1 {\e_k} V_k, \qquad
\Omega_k := \frac 1 \e_k \Omega,
\qquad
\varphi_k(x,\xi) := \varphi(\e_k x, \xi).
\]
Then by rescaling we find that
\[
Z_k \mbox{ minimizes }P_{\varphi_k}(\cdot, \R^n) \ \mbox{ in } \ 
\{ W\subset\Omega : W \setminus B(1)  = \Omega_k \setminus B(1) \} .
\]
It follows that, if $ W \setminus B(1)  = \Omega_k \setminus B(1)$, then
$P_{\varphi_k}(Z_k , U) \le P_{\varphi_k}(W;U)$ for any open set $U$ such
that $B(1)\subset\subset U$.

It is convenient to write points in $\R^n$ in the form
$(x', x_n)$ with $x'\in \R^{n-1}$. We will also write $B'(r)$ to denote the open ball of radius $r$  about the origin in $\R^{n-1}$.
Then the  smoothness of $\partial \Omega$ and our choice of coordinates
imply that there exists a smooth function $\omega: B'(\delta)\to \R$
for some $\delta>0$, such that near $x_0=0$, we can write $\partial \Omega$ 
as the graph of $\omega$, and that $\nabla \omega (0)=0$. 
If we define 
$\omega_k(x') = \frac 1 {\e_k}\omega(\e_k x')$, then for $k$ sufficiently large,
\[
\Omega_k \cap \left(B'(2)\times (-2,2) \right) = \{ (x', x_n) : x' \in B'(2),  \  \ -2<  x_n < \omega_k(x')  \}.
\]
Note that $\omega_k \to 0$ in $C^2(B'(2))$ as $k\to \infty$.
We will also write
\[
A_k ':= \{ x'\in \R^{n-1} : (x', \omega_k(x')) \in B(1) \},
\qquad
A_k := \{ x = (x', x_n)\in B_1 : x'\in A_k' \}. %A_k' \times (-1,1).
\]
%\blue{include figure?}
Then $A_k'$ is a subset of $B'(1)$ with $C^2$ boundary, and $A_k'\to B'(1)$ in the Hausdorff sense as $k\to \infty$.

%\[
%\partial \Omega_\e \cap (B'(2)\times [-2,2]) = \{ (x', \omega_\e(x')) : x' \in B'(2)  \},\qquad
%\omega_\e(x') = \omega(\e x').
%\]

\medskip

We now claim that %We will later show that 
if $k$ is large enough, then
there exists a function $v_k \in C^{1,\alpha}(\bar A'_k)$ such that 
$v_k\le \omega_k$ in $\Omega$, $v_k = \omega_k$ on $\partial \Omega$, and
\begin{equation}
%\begin{multline}
%\begin{aligned}
\mbox{ $v_k \le \omega_k$ in $A_k'$, \qquad \ \ $v_k = \omega_k$ on $\partial A_k'$,  } \ \ \qquad
%\label{vk.bc}\end{equation}
%\begin{equation}
% &\hspace{16em} 
\ Z_k \cap A_k  = \{ (x', x_n) \in A_k \ : \   \  x_n < v_k(x') \  \}.
%\end{aligned}
\label{obstacle.smooth} \end{equation}
%\end{multline} 
Toward this end, we first note that standard compactness results imply that there exists some $Z\subset \R^n$ such that after passing to a subsequence if necessary,  $\chi_{Z_k} \to \chi_Z$ in $L^1_{loc}$ as $k\to \infty$. Then standard lowersemicontinuity results and the optimality of $Z_k$
imply that, 
if we write  $\varphi_0(\xi) := \varphi(0,\xi)$ and $H^- :=\{ (x', x_n) : x_n<0\}$, 
then
\begin{align*}
P_{\varphi_0}(Z; B(2)) \le \liminf_{k\to \infty} P_{\varphi_k}(\Omega_k; B(2)) 
&\le 
 \liminf_{k\to \infty} P_{\varphi_k}(Z_k; B(2)) \\
& = P_{\varphi_0}( H^- , B(2)).
\end{align*}
Also, since $Z_k = \Omega_k$ outside $B(1)$, it is clear that $Z\setminus B(1) = H_-\setminus B(1)$.
However, convexity properties of $\varphi_0$, see $C_3$, imply that a half-plane
is always strictly $\varphi_0$-area minimizing with respect to compactly supported perturbations,
so it follows that  $Z = H^-$.
%It follows from $C_3$ that $\varphi_0$ is strictly convex in the sense of  \eqref{sconvex}, so we infer from Lemma \ref{L.convex.min} that

Once this is known, our claim about the existence of functions $v_k$, for $k$ large enough,
satisfying  \eqref{obstacle.smooth} follows from \cite{DS, White}.
More precisely, the proof of \cite[Theorem 1]{White} shows that the regularity theory for
almost-minimizing currents, together with the fact that $\chi_{Z_k} \to \chi_{H^-}$ in $L^1_{loc}$,
implies the claim. And the specific regularity results needed in our setting (where $\varphi_k$ and the boundary data for $Z_k$ depend on
$k$, but are uniformly bounded in suitable norms as $k\to \infty$)
are established in \cite[Theorem 6.1]{DS}.

For $w\in C^\infty(A_k')$ such that $w = \omega_k$ on $\partial A_k'$
and $\mbox{graph}(w) \subset A_k\subset B(1)$, we define
\[
I_k[w]:= \int_{A'_k} \varphi_k(x', w, -Dw , 1) dx'.
\]
Then  $I_k[w] = P_{\varphi_k}\left( \left\{  (x',x_n)\in A_k : x_n < w(x') \right\} \right)$ for
such $w$.
It then follows from the optimality of $Z_k$ that
\[
I_k[v_k] \le I_k[v_k - t w] \qquad\mbox{ for $w\in C^\infty(A'_k)_c$ such that 
$w(x)\ge 0$ everywhere}
\]
if $t$ is positive and sufficiently small. It follows that if $w\ge 0$ with support in $A_k'$, then
\[
0 \ge \lim_{t\searrow 0}
\frac 1t(I_k[v_k]  -  I_k[v_k - t w])
=
\int_{A'_k} \varphi_{k, x_n}(x', v_k, -Dv_k , 1)w
-
\sum_{i=1}^{n-1}\varphi_{k, \xi_i}(x', v_k, -Dv_k , 1) w_{x_i}
dx'.
\]
Thus $v_k$ satisfies
\begin{equation}
-\sum_{i=1}^{n-1} \partial_{x_i} \varphi_{k,\xi_i}(x', v_k, -Dv_k,1) -
 \varphi_{k, x_n}(x', v_k, -Dv_k , 1) \ge  0
\label{vk.subsol}\end{equation}
weakly.
On the other hand,  \eqref{barrier.suff} implies that
\begin{equation}
-\sum_{i=1}^{n-1} \partial_{x_i} \varphi_{k,\xi_i}(x', \omega_k, -D\omega_k,1) -
 \varphi_{k, x_n}(x', \omega_k, -D\omega_k , 1) \le  0,
 \label{omegak.supersol}\end{equation}
with strict inequality on a dense subset of $A_k'$.  (We recall the proof of this below, for the convenience
of the reader.)
If we let $w =v_k - \omega_k$, then arguing as
in the proof of Lemma \ref{HopfMax}, 
we find  that $w$ is a weak subsolution a linear elliptic problem
of the form \eqref{linear.diff}, i.e. that
\[
\sum_{i,j=1}^{n-1} \frac{\partial}{\partial x_i} \left(-\langle \varphi_{k,\xi_i \xi_j}\rangle w_{x_j} + \langle \varphi_{k,\xi_i x_n}\rangle w \right)
+\langle\varphi_{k,x_n \xi_i} \rangle w_{x_i} 
- \langle\varphi_{k,x_n x_n}\rangle w \ \le  \ 0 , 
\]
weakly in $A_k'$. Moreover,  \eqref{obstacle.smooth}
implies that $w\le 0$ in $A_k'$, and $w=0$ on $\partial A_k'$.

Our assumption \eqref{barrier.fail} implies that $w=0$
at some point $y_k \in A_k'$, and the weak Harnack inequality for
subsolutions then implies that $w=0$ everywhere in a small
ball about $y_k$. This implies that $v_k = \omega_k$ in this
small ball, which is impossible, in view of \eqref{vk.subsol}, 
together with the fact that \eqref{omegak.supersol} holds with
strict inequality on a dense subset. This contradiction
completes the proof of the lemma.

\hfill $\Box$

Finally we prove that \eqref{barrier.suff} and \eqref{omegak.supersol} are equivalent. To show this, it suffices to prove that if $w$ is a $C^2$ function $\R^{n-1}\to \R$ 
and if $d:\R^n\to \R$ is the signed distance to the graph of $w$ (positive below the 
graph
and negative above) then
\begin{equation}
\sum_{i=1}^n \partial_{x_i}\varphi_{\xi_i}(x, Dd)
=
\sum_{a=1}^{n-1} \partial_{x_a}\varphi_{\xi_a}(x', w, -Dw, 1)
+\varphi_{x_n}(x', w, -Dw, 1).
\label{reformulate}\end{equation}
at points $x = (x', w(x'))$ in the graph of $w$.

First, since $\varphi(x,\lambda\xi) = \lambda \varphi(x,\xi)$ for all $x$,
it follows that $\varphi_{x_n}(x,\lambda\xi) = \lambda \varphi_{x_n}(x,\xi)$, and hence that
\[
(-Dw, 1) \cdot \nabla_\xi \varphi_{x_n}(x', w, -Dw, 1) =  \varphi_{x_n}(x', w, -Dw, 1).
\]
Using this we see that
\begin{multline*}
\sum_{a=1}^{n-1} \partial_{x_a}\varphi_{\xi_a}(x', w, -Dw, 1)
+\varphi_{x_n}(x', w, -Dw, 1)\\
=
\sum_{i,j=1}^{n}
\varphi_{x_i \xi_i}(x', w, -Dw, 1)
-
\sum_{a,b=1}^{n-1} 
\varphi_{\xi_a \xi_b}(x', w, -Dw, 1)w_{x_a x_b}.
\end{multline*}
Comparing with  \eqref{reformulate},
we find  that it now suffices to prove that 
\[
\sum_{i,j=1}^{n} 
\varphi_{\xi_i \xi_j}(x, Dd) d_{x_i x_j}
\ = \ 
\sum_{a,b=1}^{n-1} 
\varphi_{\xi_a \xi_b}(x', w, -Dw, 1)w_{x_a x_b}.
\]
To do this, it is helpful to define $f(x) = x_n - w(x')$.
Then the zero level sets of $d$ and of $f$ coincide, so
$Df = \frac {|Df|}{|Dd|} Dd = |Df| Dd$. Then by
homogeneity, writing $x = (x',w(x'))$,
\[
\varphi_{\xi_a \xi_b}(x', w, -Dw, 1)
=
\varphi_{\xi_a \xi_b}(x, Df)  = |Df|^{-1}
\varphi_{\xi_a \xi_b}(x, Dd) 
\]
and thus, using the form of $f$
and the fact that $\partial_{\xi_i\xi_j}(x,\nu) \nu_j = 0$
(another consequence of homogeneity) we rewrite
\[
\sum_{a,b=1}^{n-1} 
\varphi_{\xi_a \xi_b}(x', w, -Dw, 1)w_{x_a x_b}
=
\sum_{i,j=1}^{n} 
\varphi_{\xi_i \xi_j}(x, Dd) \frac{f_{x_i x_j}}{|Df|} 
= 
\sum_{i,j=1}^{n} 
\varphi_{\xi_i \xi_j}(x, Dd) \partial_{x_j}(\frac{f_{x_i }}{|Df|} ).
\]
However, $Dd = \frac {Df}{|Df|}$ at points in the graph of $w$ (ie, the zero level-set of
both $d$ and $f$) which implies that all tangential derivatives of $Dd$ and $Df/|Df|$
are equal. Then, again appealing to the fact that  $\partial_{\xi_i\xi_j}(x,\nu) \nu_j = 0$,
we conclude that
\[
\sum_{i,j=1}^{n} 
\varphi_{\xi_i \xi_j}(x, Dd) \partial_{x_j}(\frac{f_{x_i }}{|Df|} )
=
\sum_{i,j=1}^{n} 
\varphi_{\xi_i \xi_j}(x, Dd)d_{x_i x_j},
\]
completing the proof of \eqref{reformulate}. \hfill $\Box$
\vspace{0.5cm}

\section{Acknowledgements} We are grateful to Ben Stephens for many helpful discussions. The first author was partially supported by an NSERC Discovery Grant. The second author was partially supported by MITACS and NSERC postdoctoral fellowships. The third author was partially supported by an NSERC Discovery Grant; he wishes to thank the Mittag-Leffler Institute for their wonderful hospitality.

 \end{document}